\documentstyle[amssymb,12pt]{amsart}

\newtheorem{prop}{Proposition}[section]
\newtheorem{defin}{Definition}[section]
\newtheorem{lemma}{Lemma}[section]
\newtheorem{thm}{Theorem}[section]
\newtheorem{corollary}{Corollary}[section]

\theoremstyle{remark}
\newtheorem{remark}{Remark}

\begin{document}
\newcommand{\nc}{\newcommand} \nc{\on}{\operatorname}
\nc{\pa}{\partial}
\nc{\cA}{{\cal A}}\nc{\cB}{{\cal B}}\nc{\cC}{{\cal C}}
\nc{\cE}{{\cal E}}\nc{\cG}{{\cal G}}\nc{\cH}{{\cal H}}
\nc{\cX}{{\cal X}}\nc{\cR}{{\cal R}}\nc{\cL}{{\cal L}}
\nc{\sh}{\on{sh}}\nc{\Id}{\on{Id}}\nc{\Diff}{\on{Diff}}
\nc{\ad}{\on{ad}}\nc{\Der}{\on{Der}}\nc{\End}{\on{End}}\nc{\res}{\on{res}}
\nc{\Imm}{\on{Im}}\nc{\limm}{\on{lim}}\nc{\Ad}{\on{Ad}}
\nc{\Hol}{\on{Hol}}\nc{\Det}{\on{Det}}
\nc{\de}{\delta}\nc{\si}{\sigma}\nc{\ve}{\varepsilon}
\nc{\al}{\alpha}
\nc{\CC}{{\Bbb C}}\nc{\ZZ}{{\Bbb Z}}\nc{\NN}{{\Bbb N}}
\nc{\AAA}{{\Bbb A}}\nc{\cO}{{\cal O}} \nc{\cF}{{\cal F}}
\nc{\la}{{\lambda}}\nc{\G}{{\frak g}}\nc{\A}{{\frak a}}
\nc{\HH}{{\frak h}}
\nc{\N}{{\frak n}}\nc{\B}{{\frak b}}
\nc{\La}{\Lambda}
\nc{\g}{\gamma}\nc{\eps}{\epsilon}\nc{\wt}{\widetilde}
\nc{\wh}{\widehat}
\nc{\bn}{\begin{equation}}\nc{\en}{\end{equation}}
\nc{\SL}{{\frak{sl}}}

%
%
%

\newcommand{\ldar}[1]{\begin{picture}(10,50)(-5,-25)
\put(0,25){\vector(0,-1){50}}
\put(5,0){\mbox{$#1$}} 
\end{picture}}

\newcommand{\lrar}[1]{\begin{picture}(50,10)(-25,-5)
\put(-25,0){\vector(1,0){50}}
\put(0,5){\makebox(0,0)[b]{\mbox{$#1$}}}
\end{picture}}

\newcommand{\luar}[1]{\begin{picture}(10,50)(-5,-25)
\put(0,-25){\vector(0,1){50}}
\put(5,0){\mbox{$#1$}}
\end{picture}}

\title{Elliptic quantum groups $E_{\tau,\eta}({\frak{sl}}_{2})$ 
and quasi-Hopf algebras}

\author{B. Enriquez}\address{Centre de Math\'ematiques, URA 169 du CNRS,
Ecole Polytechnique,
91128 Palai-seau, France}

\author{G. Felder}
\address{D-Math, ETH-Zentrum, HG G46, CH-8092 Z\"urich, Suisse}

\date{March 1997}
\maketitle

\begin{abstract} 
We construct an algebra morphism from the elliptic quantum group
$E_{\tau,\eta}(\frak{sl}_{2})$ to a certain elliptic version of the 
``quantum loop groups in higher genus'' studied by V. Rubtsov and the 
first author. This provides an embedding of $E_{\tau,\eta}(\frak{sl}_{2})$ 
in an algebra ``with central extension''. In particular we construct 
$L^\pm$-operators obeying a dynamical
version of the Reshetikhin--Semenov-Tian-Shansky relations.
To do that, we construct the factorization
of a certain twist of the quantum loop algebra, that automatically satisfies
the ``twisted cocycle eqaution'' of O. Babelon, D. Bernard and
E. Billey, and therefore provides a solution of the dynamical
Yang-Baxter equation. 
\end{abstract}

\subsection*{Introduction}

The aim of this paper is to compare the $\SL_{2}$-version of the
elliptic quantum groups introduced by the second author (\cite{Ell-QG})
with quasi-Hopf algebras introduced by V. Rubtsov and the
first author (\cite{HGQG,Enr-Rub}). Elliptic quantum groups are
presented by exchange (or ``$RLL$'') relations, whereas the algebras
of \cite{HGQG} are ``quantum loop algebras''. 
Our result can be viewed as an elliptic version of the results
of J. Ding and I. Frenkel (\cite{DF}) and of S. Khoroshkin (\cite{Kh}), 
where Drinfeld's quantum current algebra 
(\cite{D-new}) was shown to be isomorphic with the Reshetikhin-Semenov 
$L$-operator algebra of \cite{RS,FRT}, in the trigonometric and rational 
case respectively.  

Elliptic quantum groups are based on a matrix solution $R(z,\la)$ of the
dynamical Yang-Baxter equation (YBE). Here ``dynamical'' means that in
addition to the spectral parameter $z$ (belonging to an elliptic curve
$E$), $R$ depends on a parameter $\la$, which undergo
certain shifts in the Yang-Baxter equation. The $RLL$ relations defining the
elliptic quantum groups $E_{\tau,\eta}(\SL_{2})$ are then an algebraic
variant of the dynamical YBE. 

In \cite{BBB}, O. Babelon, D. Bernard and E. Billey studied the
relation beween the dynamical and quasi-Hopf Yang-Baxter
equations. They showed that given a family of twists of a
quasi-triangular Hopf algebra, satisfying a certain ``twisted cocycle
equation'', the quasi-Hopf YBE satisfied by the twisted $R$-matrices was
indeed equivalent to the dynamical YBE. 

The quantum loop algebras of \cite{HGQG} are generally associated with
complex curves and rational differentials. As it was shown in
\cite{Enr-Rub}, they can be endowed with a quasi-Hopf structure,
quantizing Drinfeld's ``higher genus Manin pairs'' (\cite{D-qH}). 
To precise which quantum loop algebra should be associated with elliptic
quantum groups, we first make a quasi-classical study (section
\ref{classical}). The classical
$r$-matrix $r_{\la}(z,w)$ 
associated with $R(z-w,\la)$ corresponds to what we may call a
``dynamical
Manin triple'', that is to a family $\G_{\la}$ of maximal isotropic
complements of a fixed maximal isotropic subalgebra $\G_{\cO}$ in a
Lie algebra $\G$, endowed with a nondegenerate inner product. 
Here $\G$ is a double extension of the Lie algebra $\SL_{2} \otimes
k$, where $k$ is the local field at the origin of the elliptic curve $E$
of modulus $\tau$,
$\G_{\cO}$ is a cocentral extension of $\SL_{2} \otimes \cO$, where
$\cO$ is the local ring of $E$ at the same point, and $\G_{\la}$ is an
extension of the sum $(\N_{+}\otimes L_{\la}) \oplus (\HH \otimes L_{0})
\oplus (\N_{-}\otimes L_{-\la})$, where $L_{\la}$ are the sets of
expansions at the origin of $E$ of functions on its universal cover with
certain transformation properties. 

According to \cite{Enr-Rub}, this Manin pair $(\G,\G_{\cO})$ defines
quantum loop algebras $U_{\hbar}\G_{\cO} \subset U_{\hbar}\G(\tau)$;
$U_{\hbar}\G(\tau)$ is endowed with coproducts $\Delta$ and $\bar\Delta$,
which are conjugated by a certain twist $F$. These algebras are
presented in section \ref{pres} (analogous relations can be found in
\cite{DI}). Our aim is to find a solution of the dynamical YBE in this
algebra, quantizing $r_{\la}(z,w)$. 

To do that, we will construct twisted cocycles (in the sense of
\cite{BBB}). For 
that, we follow the method of \cite{Rat}. In that paper, we gave the
construction of a Hopf algebra cocycle in the double Yangian algebra
$DY(\SL_{2})$, by factorizing the Yangian analogue $F_{Yg}$ of $F$ as a
product $F_{Yg} = F_{2}F_{1}$, with $F_{1}\in A^{<0} \otimes A^{\ge 0}$ and 
$F_{2}\in A^{\ge 0} \otimes A^{<0}$, with $A^{\ge 0}$ and $A^{<0}$ the
subalgebras of $DY(\SL_{2})$ generated by the nonnegative and negative 
Fourier modes of the quantum currents. 

This method does not apply directly here: we look for a family of
twists $(F^{1}_{\la},F^{2}_{\la})$. Moreover, we indeed have an analogue
of $A^{\ge 0}$, that is $U_{\hbar}\G_{\cO}$, but no analogue of
$A^{<0}$. Our idea is to construct subalgebras of the 
algebra of families of elements of 
$U_{\hbar}\G(\tau)^{\otimes 2}$ depending on $\la$, 
which play the role of $A^{<0}
\otimes A^{\ge 0}$ and $A^{\ge 0} \otimes A^{<0}$, and to which
$F^{1}_{\la}$ and $F^{2}_{\la}$ should belong (section \ref{A+-}). 
The properties of these
algebras $A^{-+}$ and $A^{+-}$ are based on properties of
relations between ``half-currents'' (generating series for elements of
deformations of $\N_{\pm}\otimes L_{\pm\la}$ and $\N_{\pm}\otimes \cO$),
following from the vertex relations
for $U_{\hbar}\G(\tau)$ (sections \ref{ee}, \ref{ee:sh}). 

The decomposition $F$ as a product $F^{2}_{\la}F^{1}_{\la}$ is carried
out in section \ref{sect-decomp}. It imitates the similar decomposition
in \cite{Rat}: applying certain projections to the decomposition
identity leads us to guess the values of $F^{1,2}_{\la}$. We then show
that the decomposition identity is indeed satisfied. The proof uses
Hopf algebra duality results of section \ref{dualty}, and results on 
coproducts of section \ref{copdt}. 

Next we prove the the $F^{1}_{\la}$ obtained that way indeed satisfies
a twisted cocycle equation. For that, we introduce subalgebras of
$U_{\hbar}\G(\tau)^{\otimes 3}$, analogous to 
$A^{\ge 0, <0} \otimes DY(\SL_{2})^{\otimes 3}$
and $DY(\SL_{2})^{\otimes 3} \otimes A^{\ge 0,<0}$, and show that they
contain images by $\Delta$ and $\bar\Delta$ of the $A^{+-},
A^{-+}$ (Props. \ref{egill}, \ref{snorri}). 
This shows that the ratio $\Phi_{\la}$ between the two
sides of the twisted cocycle equation has a special form ($\sum_{i}
1\otimes a_{i} \otimes h^{i}$, where $h$ is the ``Cartan element'' of
$U_{\hbar}\G(\tau)$). We then prove a ``twisted pentagon equation'' for
$\Phi_{\la}$ (prop. \ref{compat}), which in fact shows that it is equal
to $1$. 

After that we can apply the result of \cite{BBB}, and obtain a solution
$\cR_{\la}$
of the dynamical YBE in $U_{\hbar}\G(\tau)^{\otimes 2}$
(Thm. \ref{Thm:DYBE}). We next study
level zero representations of $U_{\hbar}\G(\tau)$; this study was led in
the general case in \cite{HGQG}. We obtain a family $\pi_{\zeta}$ of
$2$-dimensional evaluation modules (Prop. \ref{fd:repres}), 
indexed by a point $\zeta$ of the
formal neighborhood of the origin in $E$; we compute the image of
$\cR_{\la}$ by these representations, and find an answer closely
connected to $R(z,\la)$. This enables us to prove that the $L$-operators
$(\pi_{\zeta} \otimes 1)(\cR_{\la})$ and $(1\otimes
\pi_{\zeta})(\cR_{\la})$ satisfy the following dynamical version of the
Reshetikhin--Semenov-Tian-Shansky relations (Thm. \ref{thm:RLL}): 
$$
R^{\pm}(\zeta-\zeta',\la) L^{\pm(1)}_{\la-\gamma h^{(2)}} (\zeta)
L^{\pm(2)}_{\la} (\zeta')
=
L^{\pm(2)}_{\la-\gamma h^{(1)}} (\zeta') L^{\pm(1)}_{\la} (\zeta)
R^{\pm}(\zeta-\zeta',\la-\gamma h) 
$$
\begin{align*}
L^{-(1)}_{\la}(\zeta) & R^{-}(\zeta-\zeta',\la-\gamma h)
L^{+(2)}_{\la}(\zeta')
\\ & \nonumber =
L^{+(2)}_{\la-\gamma h^{(1)}} (\zeta') R^{-}(\zeta-\zeta'+K\gamma,\la)
L^{-(1)}_{\la-\gamma h^{(2)}}(\zeta) {{A(\zeta,\zeta'-K\gamma)}
\over{A(\zeta,\zeta')}} , 
\end{align*}
where $R^{\pm}(z,\la)$ are 
elliptic $R$-matrices and $A$ is a certain elliptic version of the usual
ratio of gamma-functions. 
These relations extend the $RLL$-relations
of the elliptic quantum group, which we recover in
Thm. \ref{comparison}. 

Let us now say some words about problems connected to the present
work. In \cite{Ell-QG,F-Varch}, quantum Knizhnik-Zamolodchikov-Bernard
equations were defined as difference equations involving the dynamical
$R$-matrix $R(z,\la)$. It would be interesting to derive these
equations from considerations involving coinvariants of
$U_{\hbar}\G(\tau)$-modules. This may
also shed light on the question what the equation for dependence in the
moduli should be in the quantum situation. There is some indication
that this equation is the Ruijsenaars-Schneiders (RS) equation at the
critical level. At that level, one
may apply the Reshetikhin-Semenov method for expressing elements of the
center of $U_{\hbar}\G(\tau)$, and then explicitly compute their
actions on coinvariant spaces. Recent work of A. Varchenko and the second
author leads to the impression that the situation is more complicated
outside the critical level. It might be interesting to connect such an
approach to the RS models with those of \cite{ACF,ABB}. 

Finally, it would be interesting to find an analogue of the theory
developed in the present work, for the situations of higher genus
(\cite{HGQG}). Recall that the classical $r$-matrices underlying the
elliptic quantum groups are dynamical $r$-matrices for the Hitchin
system associated with an elliptic curves. Dynamical $r$-matrices for
Hitchin systems in higher genus have been introduced in
\cite{Houches,Hab}. In this respect, it seems that a dynamical version
of the Poisson-Lie theory would be of interest. 

Let us also mention here the work \cite{Fr}, where a dynamical approach
to other `` elliptic quantum groups'' (those of \cite{Fo}) is
presented. 

This work was done during our stay at the ``Semestre
syst\`emes int\'egrables'' organized at the Centre
Emile Borel, Paris, UMS 839, CNRS/UPMC. 
We would like to express our thanks to its organizers for their
invitation to this very stimulating meeting. We also would like to
acknowledge discussions with O. Babelon, D. Bernard, C. Fronsdal,
V. Rubtsov and A. Varchenko on the subject of this work. 

\section{Quantum loop algebras associated with elliptic curves} 

\subsection{The classical situation} \label{classical}

In \cite{Enr-Rub}, we constructed quasi-Hopf algebras, associated to the
general data of a Frobenius algebra, a maximal isotropic
subalgebra of it, and an invariant derivation. 
An example of such data is the following. 

Let us fix a complex number $\tau$, with $\Imm(\tau)>0$. Let $L\subset\CC$
be the lattice $\ZZ + \tau\ZZ$; call $E$ the elliptic curve $\CC/L$. Let
$z$ be the coordinate on $\CC$, and let $\omega$ be the differential
form on $E$, equal to $dz$. 
Let $k= \CC((z))$ 
be the completed local field of $E$ at its origin $0$, and $\cO =
\CC[[z]]$ 
the completed local ring at the same point. Endow $k$ with the scalar
product $\langle , \rangle_{k}$ defined by 
$$
\langle f,g \rangle_{k} = \res_{0}(fg\omega). 
$$
Define on $k$ the derivation $\pa$ to be equal to $d/dz$. Then $\pa$ is
invariant w.r.t $\langle, \rangle_{k}$, and
$\cO$ is a maximal isotropic subring of $k$. 

Let us set $\A = \frak{sl}_{2}(\CC)$, and denote by $\langle ,
\rangle_{\A}$ an invariant scalar product on 
$\A$.  Let us set $\G = (\A \otimes
k) \oplus \CC D \oplus \CC K$; let us define on $\G$ the Lie algebra
stucture defined by the central extension of $\A \otimes k$
$$
c(x\otimes f,y\otimes g) = \langle x,y \rangle_{\A} \langle f,\pa g
\rangle_{k} K
$$
and by the derivation $[D,x\otimes f] = x \otimes \pa f$. 

Let $\G_{\cO}$ be the Lie subalgebra of $\G$ equal to 
$(\A \otimes \cO) \oplus \CC D$.  
Define $\langle , \rangle_{\A\otimes k}$ as the tensor product
of $\langle , \rangle_{\A}$ and $\langle , \rangle_{k}$, and $\langle ,
\rangle_{\G}$ as the scalar product on $\G$ defined by 
$\langle , \rangle_{\G}|_{\A \otimes k} = \langle , \rangle_{\A \otimes
k}$, $\langle D, \A \otimes k \rangle_{\G} = \langle K, \A \otimes k
\rangle_{\G} = 0$, and $\langle D,K\rangle_{\G}=1$. 
Then $\G_{\cO}$ is a maximal isotropic Lie subalgebra of $\G$. 

To define maximal isotropic supplementaries of $\G_{\cO}$ in $\G$, 
we first define certain subspaces of $k$. 

For $\la\in\CC$, define $L_{\la}$ as follows. If $\la$ does not belong to 
$L$, define $L_{\la}$ to be the set of expansions 
near $0$ of all holomorphic functions on $\CC - L$, 
$1$-periodic and such that $f(z + \tau) = e^{ - 2 i \pi \la} f(z)$. 
For $\la =0$, 
define $L_{\la}$ as the maximal isotropic subspace of $k$ containing all
holomorphic functions $f$ on $\CC - L$, 
$L$-periodic, such that $\oint_{a} f(z) dz =0$, where $a$
is the cycle $(i\eps,i\eps+1)$ (with $\eps$ small and $>0$). 
Finally, define $L_{\la} = e^{-2i\pi m z}L_{0}$ for $\la
= n + m\tau$.  

Let $\theta$ be the holomorphic function defined on $\CC$ by the 
conditions that $\theta'(0) = 1$, 
the only zeroes of $\theta$ are the points of $L$,  
$\theta(z+1) = - \theta(z)$, and
$\theta ( z + \tau) = - e^{ - i\pi\tau} e^{ - 2 i \pi z }\theta(z)$. 
$\theta$ is then odd. 

We then have 
\begin{equation} \label{L}
L_{\la} = \oplus_{j\ge 0} \CC \left( { {\theta ' }\over
\theta} \right)^{(j)} e^{ - 2i\pi m z}, \quad \on{if} \quad \la = n + m
\tau,  
\end{equation} 
\begin{equation} \label{Lla}
L_{\la} = \oplus_{i\ge 0}
\CC \left( {{ \theta ( \la + z )} \over
{ \theta(z)} }\right)^{(i)}, \quad \on{if} \quad \la\in \CC - L, 
\end{equation}
where for $f\in k$, we let $f' = \pa f$ and $f^{(i)} = \pa^{i}f$.  

Moreover, the orthogonal of $L_{\la}$ for the scalar product $\langle , 
\rangle_{k}$ is equal to $L_{-\la}$. 

Consider now the decomposition 
\begin{equation} \label{decomp}
\G = \G_{\cO} \oplus
\G_{\la}, 
\end{equation}
where 
\begin{equation}
\G_{\la} = 
({\frak h} \otimes L_{0}) \oplus (\N_{+}\otimes L_{\la}) \oplus 
(\N_{-}\otimes L_{-\la})
\oplus \CC K . 
\end{equation}
$\G_{\cO}$ is a maximal isotropic subalgebra of 
${\frak g}$, and $\G_{\la}$ is a maximal isotropic subspace 
of it. Therefore, (\ref{decomp}) defines a Lie quasi-bialgebra structure 
on  ${\frak g}_{\cO}$, and 
(as in \cite{Enr-Rub}), of double Lie quasi-bialgebra
on ${\frak g}$. Its classical 
$r$-matrix is given by the formula
$$
r_{\la} = D \otimes K + \sum_{i} {1\over 2}
h[e^{i}] \wedge h[e_{i;0}] 
+
e [e^{i}] \wedge f [e_{i; \la}] 
+
f [e^{i}] \wedge e[e_{i; -\la}] , 
$$
for $(e^{i})_{i\ge 0},(e_{i;\la})_{i\ge 0}$ 
dual bases of $\cO$ and $L_{\la}$, with the valuation of $e^{i}$ tending 
to infinity with $i$, and we denote $x\otimes f$ by $x[f]$; 
in other terms,
\begin{align} \label{class-r}
& r_{\la}(z,w) = 
{1\over 2}(h\otimes h){{\theta'}\over{\theta}}(z-w)
+
(e\otimes f){{\theta(z-u+\la)}\over{\theta(z-u)\theta(\la)}} 
+
(f\otimes e){{\theta(z-u-\la)}\over{\theta(z-u)\theta(-\la)}} 
\\ & \nonumber + D\otimes K ; 
\end{align}
this formula (without $D\otimes K$) 
coincides with that of the classical $r$-matrix arising 
in the elliptic versions of the KZB equations (see \cite{FW}) and of the 
Hitchin system (see \cite{ER-Hitch}). 

\begin{remark}
$\cO$ plays the role of the ring $R$, in the notation of
\cite{Enr-Rub}. 
\end{remark}

\subsection{Relations for $U_{\hbar}\G(\tau)$} \label{pres}

The quasi-Hopf algebra associated in \cite{Enr-Rub} to the Lie
quasi-bialgebras $\G$ and $\G_{\cO}$, are twists of the Hopf algebra
$(U_{\hbar}\G(\tau),\Delta)$ 
and of its subalgebra $U_{\hbar}\G_{\cO}(\tau)$, that we
now present. We will sometimes denote $U_{\hbar}\G(\tau)$ by $A(\tau)$ or
simply by $A$, and $U_{\hbar}\G_{\cO}$ by $A^{+}$. 

Let $\hbar$ be a formal parameter. Generators of $U_{\hbar}\G(\tau)$
are $D,K$ and the $x[\eps]$, $x=e,f,h$, $\eps\in k$; they are subject to
the relations 
$$
x[\al\eps] = \al x[\eps],\quad x[\eps + \eps'] = x[\eps] +
x[\eps'] , \quad \al\in\CC, \eps,\eps'\in k. 
$$
They serve to define the generating series 
$$
x(z) = \sum_{{i\in\ZZ}} x[\eps^{i}]\eps_{i}(z), \quad x = e,f,h,  
$$
$(\eps^{i})_{i\in\ZZ},(\eps_{i})_{i\in\ZZ}$ dual bases of $k$; recall
that 
$(e^{i})_{i\in\NN},(e_{i;0})_{i\in\NN}$ are dual bases of $\cO$ and $L_{0}$
and set 
$$
h^{+}(z) = \sum_{i\in\NN} h[e^{i}] e_{i;0}(z), \quad 
h^{-}(z) = \sum_{i\in\NN} h[e_{i;0}] e^{i}(z).  
$$
We will also use the series 
$$
K^{+}(z) = e^{({{q^{\pa}-q^{-\pa} }\over{2\pa}} h^{+})(z)}, \quad
K^{-}(z) = q^{h^{-}(z)}, 
$$
where $q = e^{\hbar}$. 
The relations presenting $U_{\hbar}\G(\tau)$ are then 
\begin{equation}  \label{K-K}
[K^{+}(z) , K^{+}(w)] = [K^{-}(z) , K^{-}(w)] = 0, 
\end{equation}
\begin{align} \label{K+K-}
& \theta(z-w-\hbar) \theta(z-w+\hbar + \hbar K)
K^{+}(z)K^{-}(w) 
\\ & \nonumber 
 = \theta(z-w+\hbar) \theta(z-w-\hbar+\hbar K)
K^{-}(w) K^{+}(z) , 
\end{align}
\begin{align} \label{K+:e}
K^{+}(z)  e (w) K^{+}(z)^{-1} = { {\theta(z-w+\hbar)} \over
{\theta(z-w-\hbar)} } & e(w)
\end{align}
\begin{align} \label{K-:e}
K^{-}(z) e (w) K^{-}(z)^{-1} = { {\theta(w-z +\hbar K+\hbar)} \over
{\theta(w-z+\hbar K-\hbar)} 
} e(w),
\end{align}
\begin{equation} \label{K-f}
K^{+}(z)  f (w) K^{+}(z)^{-1} = { {\theta(w-z+\hbar)} \over
{\theta(w-z-\hbar)}} f(w),
K^{-}(z) f (w) K^{+}(z)^{-1} = { {\theta(z-w +\hbar)} \over
{\theta(z-w-\hbar)}}f(w),
\end{equation}
\begin{equation} \label{vertex-e}
\theta(z-w-\hbar) e(z) e(w) =  {\theta(z-w+\hbar)} e(w)e(z) , 
\end{equation}
\begin{equation} \label{vertex-f}
\theta(w-z-\hbar) f(z) f(w) 
= \theta(w-z+\hbar) f(w) f(z) , 
\end{equation}
\begin{equation} \label{e-f-ell}
[e(z) , f(w)] = {1 \over \hbar} \left( \delta(z,w) K^{+}(z) -
\delta(z,w-\hbar K) K^{-}(w)^{-1}\right). 
\end{equation}

Let us introduce the generating series $k^{+}(z)$ and  $k^{-}(z)$,
defined by 
\begin{equation}
k^{+}(z) = e^{({{ q^{\pa} - 1}\over{2\pa}} h^{+})(z) }, 
\quad
k^{-}(z) = q^{({ 1 \over{ 1 + q^{-\pa}}} h^{-})(z) } ; 
\end{equation}
they satisfy the relations 
\begin{equation} \label{rel-k}
K^{+}(z) = k^{+}(z) k^{+}(z-\hbar), \quad
K^{-}(z) = k^{-}(z) k^{-}(z-\hbar). 
\end{equation}
Equations (\ref{K+:e}), (\ref{K-:e}) and (\ref{K-f}) may be replaced by 
\begin{align} \label{k-e}
k^{+}(z)  e (w) k^{+}(z)^{-1} = { {\theta(z-w+\hbar)} \over
{\theta(z-w)}} & e(w), 
\\ & \nonumber
k^{-}(z)  e (w) k^{-}(z)^{-1} = { {\theta(w-z + \hbar K)} \over
{\theta(w-z+ \hbar K-\hbar)}} e(w),
\end{align}
and
\begin{equation} \label{k-f}
k^{+}(z) f (w) k^{+}(z)^{-1} = { {\theta(w-z)} \over
{\theta(w-z-\hbar)}} f(w),
\quad
k^{-}(z)   f (w) k^{-}(z)^{-1} = { {\theta(z-w +\hbar)} \over
{\theta(z-w)}} f(w);
\end{equation}
moreover, we will have 
\begin{equation} \label{k-k}
(k^{\pm}(z), k^{\pm}(w)) = 1, \quad
(k^{+}(z) , k^{-}(w)) = f_{K}(z-w),
\end{equation} 
where we use the group commutator notation $(x,y) = xyx^{-1}y^{-1}$ 
and $f_{K}(\zeta)$ is the formal series of $ 1 + \hbar \CC[[\zeta]]
[[\hbar]]$ defined by the functional equation
\begin{equation} \label{funct-eq} 
f_{K}(\zeta)f_{K}(\zeta - \hbar)
=
{ {\theta{ (\zeta+\hbar)}} \over {\theta(\zeta)}} : 
{ {\theta{ (\zeta+\hbar +\hbar K)}} \over {\theta(\zeta +\hbar K)}}  ; 
\end{equation} 
if $K = - 2 p$, with $p$ integer, we have
\begin{equation} 
f_{K}(\zeta)
=
{ 
{\theta(\zeta)
\left ( \prod_{ k=1}^{p-1} \theta(\zeta - 2k\hbar)
\right)^2
\theta(\zeta-2p\hbar)
} 
\over 
{
\left ( \prod_{ k=0}^{p-1} \theta(\zeta - (2k + 1)\hbar)
\right)^2  
} } . 
\end{equation}

The algebra $U_{\hbar}\G(\tau)$ is endowed with a Hopf structure given
by the coproduct $\Delta$ defined by 
\begin{equation} \label{Delta:k:ell}
\Delta(k^{+}(z)) = k^{+}(z) \otimes k^{+}(z), \quad
\Delta(K^{-}(z)) = K^{-}(z) \otimes K^{-}(z +\hbar K_{1}), 
\end{equation}
\begin{equation} \label{Delta:e:ell}
\Delta(e(z)) = e(z)\otimes K^{+}(z) + 1\otimes e(z),
\end{equation}
\begin{equation} \label{Delta:f:ell}
\Delta(f(z)) = f(z)\otimes 1 + K^{-}(z)^{-1} \otimes f(z+\hbar K_{1}),
\end{equation}
\begin{equation} \label{Delta:D:K:ell}
\Delta(D) = D \otimes 1 + 1\otimes D, \quad \Delta(K) = K \otimes 1 + 
1 \otimes K , 
\end{equation}
the counit $\varepsilon$, and the antipode $S$ defined by them; we set
$K_{1} = K\otimes 1, K_{2}=1\otimes K$. 

$U_{\hbar}\G(\tau)$ is also endowed with another Hopf structure given
by the coproduct $\bar\Delta$ defined by
\begin{equation} \label{bar:Delta:k:ell}
\bar\Delta(k^{+}(z)) = k^{+}(z) \otimes k^{+}(z), \quad
\bar\Delta(K^{-}(z)) = K^{-}(z) \otimes K^{-}(z + \hbar K_{1}), 
\end{equation}
\begin{equation} \label{bar:Delta:e:ell}
\bar\Delta(e(z)) = e(z -\hbar K_{2})\otimes K^{-}(z -\hbar K_{2})^{-1} 
+ 1\otimes e(z),
\end{equation}
\begin{equation} \label{bar:Delta:f:ell}
\bar\Delta(f(z)) = f(z)\otimes 1 + K^{+}(z) \otimes f(z),
\end{equation}
\begin{equation} \label{bar:Delta:D:K:ell}
\bar\Delta(D) = D \otimes 1 + 1\otimes D, \quad \bar\Delta(K) = K \otimes 1 + 
1 \otimes K , 
\end{equation}
the counit $\varepsilon$, and the antipode $\bar S$ defined by them. 

The Hopf structures associated with $\Delta$ and $\bar\Delta$ 
are connected by a twist 
\begin{equation} \label{F}
F = \exp\left(\hbar \sum_{i\in\ZZ} \eps_{i} \otimes \eps^{i} \right) , 
\end{equation}
where $(\eps^{i})_{\in\ZZ}$ is the basis of $k$ dual to 
$(\eps_{i})_{i\in\ZZ}$ w.r.t. $\langle , \rangle_{k}$; that is, we have 
$\bar\Delta = \Ad(F)\circ \Delta$. 

Then $F$ satisfies the cocycle equation 
\begin{equation} \label{F-cocycle}
(F\otimes 1)(\Delta\otimes 1)(F) = (1\otimes F)(1\otimes \Delta)(F)
\end{equation}
(see \cite{Enr-Rub}, Prop. 3.1). 

\begin{remark}
In the
notation of \cite{Enr-Rub}, 8.2, we have 
$$
q(z,w) = { {\theta(z-w+\hbar)} \over {\theta(z-w-\hbar)} }, 
$$
and $\sigma(z) = z -\hbar K$. We also have the relation 
$$
[h^{+}(z), h^{-}(w)] = {1\over \hbar} \left( {\theta'\over\theta}(z-w)
- {\theta'\over\theta}(z-w+\hbar K)\right).
$$

\end{remark}

\subsection{Completions} The algebras and vector spaces introduced 
above possess natural topologies: the field $k$ and the ring $\cO$ 
are given the formal series topology; on the other hand, the
spaces $L_{\la}$ are given the discrete topology. For $V,W$ two
topological vector spaces, with basis of neighborhoods of the origin 
$(V_{a})_{a\in\ZZ}$ 
and $(W_{b})_{b\in\ZZ}$, define $V \hat\otimes W$ as the inverse limit
of the $V\otimes W / V_{a} \otimes W_{b}$, and $V\bar\otimes W$ as the inverse
limit of the $(V\otimes W)/(V_{a}\otimes W + V \otimes W_{b})$.

Define then the completed tensor algebra $T\ \hat{}(V)$ 
of $V$ as the direct sum 
$\oplus_{i\ge 0} V^{\hat\otimes i}$. Then $U_{\hbar}\G(\tau)$ is
viewed as a quotient of $T\ \hat{}(\G)$, and is endowed with the 
corresponding topology. $U_{\hbar}\G_{\cO}$ is then a closed subspace 
of $U_{\hbar}\G(\tau)$. 
On the other hand, the fields $x(z)$ belong to the completed
tensor product $U_{\hbar}\G(\tau) \bar\otimes k$. The coproduct $\Delta$
maps $U_{\hbar}\G(\tau)$ to the completion of its tensor square defined
as the suitable quotient of $T\ \hat{}(\G \oplus \G)$. 
In the sequel we will consider tensor product of subspaces of
$U_{\hbar}\G(\tau)$ to be completed w.r.t. the topology of $T\ \hat{}(\G
\oplus \G), T\ \hat{}(\G \oplus \G\oplus \G)$, etc.  

\subsection{Relations for half-currents} \label{Ug+} \label{ee}

Fix a complex number $\la$ and set for $x = e,f,$ $K^+$, 
\begin{equation} \label{x+}
x^{+}_{\la}(z) = \sum_{i} x[e^{i}]e_{i;\la}(z), 
\end{equation}
and for $x=e,f,K^{-}$, 
\begin{equation} \label{x-}
x^{-}_{\la}(z) = \sum_{i} x[ e_{i;-\la} ] e_{i}(z);  
\end{equation}
recall that $(e^{i}),(e_{i;\la})$ are dual bases of $\cO$ and $L_{\la}$. 

The fields $e(z)$ and $f(z)$ are then split according to 
\begin{equation} \label{split}
e(z) = e^{+}_{\la}(z) + e^{-}_{\la}(z), \quad
f(z) = f^{+}_{-\la}(z) + f^{-}_{-\la}(z);   
\end{equation}
we call the expression $x^{\pm}_{\la}(z)$ ``half-currents''. In the
above equality, we made use of the continuous inclusions of
$U_{\hbar}\G_{\cO} \bar\otimes L_{\la}$ and of $U_{\hbar}\G \bar\otimes
\cO$ into $U_{\hbar}\G \bar\otimes k$. 

For $x = e,f$, we have $x^{-}_{\la}(z) \in U_{\hbar}\G(\tau) 
\bar \otimes \cO$, so that 
$x^{-}(z)$ can be viewed as a formal series in $z$, regular at $0$, and 
$x^{+}_{\la}(z) \in U_{\hbar}\G(\tau) \bar\otimes L_{\la}$, 
so that $x^{+}_{\la}(z)$ can be viewed 
as a function of $z$, satisfying
\begin{equation} \label{x-per}
x^{+}_{\la}(z+1) = x^{+}_{\la}(z), 
\quad
x^{+}_{\la}(z+ \tau ) = e^{ - 2 i \pi \la}x^{+}_{\la}(z).
\end{equation}

We then have: 

\begin{prop} \label{rels}
The generating series $e^{\pm}_{\la}(z), f^{\pm}_{\la }(z)$  
satisfy the following relations: 
\begin{align} \label{eeps:eeps'}
{ {\theta(z-w-\hbar)} \over {\theta(z-w)} } & e^{\eps}_{\la+\hbar}(z)  
e^{\eps'}_{\la-\hbar}(w)
+
\eps\eps'
{{\theta(w-z-\la)\theta(-\hbar)} \over {\theta(w-z)\theta(-\la)} } 
e^{\eps'}_{\la+\hbar}(w) e^{\eps'}_{\la-\hbar}(w)
\\ \nonumber & =
{ {\theta(z-w+\hbar)} \over {\theta(z-w)} } e^{\eps'}_{\la+\hbar}(w) 
e^{\eps}_{\la-\hbar}(z)
+
\eps\eps'
{ {\theta(z-w-\la)\theta(-\hbar)} \over {\theta(z-w)\theta(-\la)} } 
e^{\eps}_{\la+\hbar}(z) e^{\eps}_{\la-\hbar}(z),
\end{align}
\begin{align} \label{feps:feps'}
{ {\theta(z-w+\hbar)} \over {\theta(z-w)} } & f^{\eps}_{\la-\hbar}(z)  
f^{\eps'}_{\la+\hbar}(w)
+
\eps\eps'{ {\theta(w-z-\la)\theta(\hbar)} \over {\theta(w-z)\theta(-\la)} } 
f^{\eps'}_{\la-\hbar}(w) f^{\eps'}_{\la+\hbar}(w)
\\ \nonumber & =
{ {\theta(z-w-\hbar)} \over {\theta(z-w)} } f^{\eps'}_{\la-\hbar}(w) 
f^{\eps}_{\la+\hbar}(z)
+
\eps\eps'
{ {\theta(z-w-\la)\theta(\hbar)} \over {\theta(z-w)\theta(-\la)} } 
f^{\eps}_{\la-\hbar}(z) f^{\eps}_{\la+\hbar}(z) ,  
\end{align}
where $\eps,\eps'$ take the values $+,-$. 
In these relations, the expressions of the form 
${{1}\over{z-w}}(f'(z,w)x^{+}_{\la}(z) - f''(z,w)x^{+}_{\la}(w))$, resp.
${{1}\over{z-w}}(f'(z,w)x^{-}_{\la}(z) - f''(z,w)x^{-}_{\la}(w))$, 
where
$f',f''$ are formal series in $z,w$ coinciding for $z=w$, and
$x=e,f,K^{+}$, should be understood as the sums 
$$
\sum_{i\ge 0} x[e^{i}] {{f'(z,w)e_{i,\la}(z) -
f''(z,w)e_{i,\la}(w)}\over{z-w}}, 
\on{resp.} 
\sum_{i<0} x[e_{i,\la}] {{f'(z,w)e^{i}(z) -
f''(z,w)e^{i}(w)}\over{z-w}}, 
$$ 
which belong to $U_{\hbar}\G(\tau) \otimes (k\otimes k)$. 
\end{prop}

{\em Proof. \/}
Let us show relation (\ref{eeps:eeps'}) in the case $\eps=\eps'=+$. Let 
us denote by $Z_{++}(z,w)$ the difference of the left and right hand 
sides of this equation, and by $\ell$ any continuous linear form 
on $A(\tau)$. Clearly, $\ell(Z_{++}(z,w))$ belongs to $L_{\la} \otimes 
L_{\la}$ (recall that $\ell[e_{i}e_{j}]$ is equal to zero when $i$ or
$j$ are large enough), and is antisymmetric.
(We attach coordinates $z,w$ to the first and second factor of the tensor 
product.)
On the other hand, the difference of $\ell( \theta(z-w)Z_{++}(z,w) ) $ and of
$\ell( \theta(z-w-\hbar)e(z)e(w)- \theta(z-w+\hbar)e(w)e(z))$ 
can be expressed as a sum of quadratic monomials in the $e^{\pm}_{\la\pm
\hbar}(z,w)$ using at least one $e^{-}_{\la\pm \hbar}(z,w)$, and
therefore belongs to $k \bar \otimes \cO + \cO \bar \otimes k$. 
Therefore, the same is true for $(z-w)\ell(Z_{++}(z,w))$. 
Let us set $(z-w)Z_{++}(z,w) = Y_{1}+ Y_{2}$, with $Y_{1} \in k \bar\otimes 
\cO, Y_{2} \in \cO \bar\otimes k$. 
Let us denote by a tilde the exchange of arguments $z$ and $w$, and set 
$Y = (Y_{1}+\wt Y_{2})/2$. We have $Y \in k \bar\otimes \cO$ and
$$
(z-w)\ell( Z_{++}(z,w)) = Y +\wt Y. 
$$  

Set now $Z = Y \sum_{i\ge 0} z^{-i-1}w^{i}$; since $Y$ belongs to 
$k \bar\otimes \cO = k[[w]]$, so does $Z$. Since we have 
$$
(z-w) [\ell(Z_{++}(z,w)) - Z + \wt Z] = 0,  
$$
it follows that for some $f$ in $k$, we have  
\begin{equation} \label{myra}
\ell(Z_{++}(z,w)) - Z + \wt Z = f(z) \delta(z,w) .   
\end{equation}
Consider both sides of (\ref{myra}) as the kernel of some operator, 
defined by $(Tf_{0})(z) = \res_{0} (f_{0}(z)\delta(z,w)dw)$. 
Since the l.h.s. of (\ref{myra}) is antisymmetric in $z$ and $w$, this 
operator should be antisymmetric, i.e. satisfy $\langle Tf_{0},g_{0} 
\rangle_{k}
+ \langle f_{0}, Tg_{0}\rangle_{k} =0$, for $f_{0},g_{0}\in k$. 
On the other hand, $T$
coincides with the multiplication operator by $f$. It follows that 
$f=0$.  Therefore $\ell(Z_{++}(z,w)) = Z - \wt Z$. Since the left and
right hand sides of this equality belongs to $L_{\la} \otimes L_{\la}$ 
and to $\cO \bar\otimes k + k \bar\otimes \cO$ respectively, 
$\ell(Z_{++}(z,w)) = 0$. This proves (\ref{eeps:eeps'}) in the case where 
$\eps = \eps' =+$. 

Let us now show (\ref{eeps:eeps'}) in the case where 
$\eps = \eps' =-$. 
Let us denote by $Z_{--}(z,w)$ the difference of the left and right hand 
sides of this equation.
Clearly, $\ell(Z_{--}(z,w))$ belongs to $(\cO \hat\otimes \cO
)[[\hbar]]$, and is antisymmetric. On the other hand, the difference 
of $\theta(z-w) \ell(Z_{--}(z,w))$ with 
$\ell( \theta(z-w-\hbar)e(z)e(w)- \theta(z-w+\hbar)e(w)e(z))$ 
belongs to $\cF_{\la, *} + \cF_{*, \la}$, where $\cF_{\la,*}$ is the subspace
of $\Hol(\CC - L)((w))$ formed by the functions $f(z,w)$ such that 
$$
f(z+1,w) = f(z,w), \quad f(z+\tau,w) = -e^{-i\pi\tau}e^{-2i\pi\la}e^{-2i\pi
(z-w)}f(z,w),
$$
any 
$\cF_{*,\la}$ is the subspace
of $\Hol(\CC - L)((z))$ formed by the functions $g(z,w)$ such that 
$$
g(z,w+1) = g(z,w), \quad g(z,w+\tau) = -e^{-i\pi\tau}e^{-2i\pi\la}e^{-2i\pi
(w-z)}g(z,w);
$$ 
here $\Hol(\CC - L)$ is the space of holomorphic 
functions defined on $\CC - L$. 
Set $\theta(z-w) \ell(Z_{--}(z,w)) = Y'_{1} + Y'_{2}$, with $Y'_{1}\in 
\cF_{\la,*}$ and $Y_{2} \in \cF'_{*,\la}$. Set $Y'=(Y'_{1}+\wt
Y'_{2})/2$. We have $Y' \in \cF_{\la,*}$ and 
$$
\theta(z-w)\ell( Z_{++}(z,w)) = Y' +\wt Y'.  
$$  
Set now 
$$
Z' = Y' \sum_{i\ge 0}(\theta^{-1})^{(i)}(z){{(-w)^{i}}\over{i!}} ;
$$ 
then $Z'$ belongs to $\Hol(\CC - L)((w))$, and we have as before 
$$
\ell(Z_{--}(z,w)) = Z' - \wt Z' + f'(z) \delta(z,w), 
$$
for some $f'\in k$. The same reasoning as above shows that $f'=0$ and then 
that $\ell(Z_{--}(z,w)) = 0$. 

Let us now prove (\ref{eeps:eeps'}) in the case where $\eps = +, \eps'=-$. 
Let us denote by $Z_{+-}(z,w)$ the difference of the left and right 
hand sides of this equality. 
Let us substract from (\ref{vertex-e}), the sum of equalities 
(\ref{eeps:eeps'}) with $\eps = \eps' = +$ and $\eps = \eps' = -$. 
We obtain that 
$Z_{+-}(z,w) + Z_{+-}(w,z) = 0$. 
Therefore $Z_{+-}(z,w)$ is antisymmetric in $z$ and $w$. 
On the other hand, we have for any linear functional $\ell$ on $A(\tau)$,
$\ell (Z_{+-}(z,w))
\in L_{\la} \bar\otimes \cO$.
Since the intersection of $L_{\la}$ and $\cO$ is zero, $Z_{+-}(z,w)$ is
equal to zero. Therefore (\ref{eeps:eeps'}) is valid in the case
$\eps = +, \eps'=-$. 

The case $\eps = -, \eps'=+$ is obtained from $\eps = +, \eps'=-$, 
by exchanging $z$ and $w$. 

Relations (\ref{feps:feps'}) can be obtained in a similar way. 
\hfill \qed \medskip

Let us define $U_{\hbar}\N_{+}(\tau)$  and $U_{\hbar}\N_{-}(\tau)$ 
as the subalgebras of $A(\tau)$ generated by the $e[\eps],\eps\in k$ 
and the $f[\eps],\eps\in k$. 

Let us denote by $U^{(e)}_{+}(\tau)$, resp.
$U^{(f)}_{+}(\tau)$ the subalgebras of 
$U_{\hbar}\N_{\pm}(\tau)$ generated by 
the $e[r],r\in \cO$, resp. $f[r],r\in\cO$. 

For $\mu\in \CC - L$, and $x = e,f,h$, denote by 
$x^{-}_{\mu}[\eps]$ the element 
$x[p_{-\mu}(\eps)]$, where 
$p_{-\mu}$ is the projection on $L_{-\mu}$ parallel to $\cO$.
For $\beta\in \CC$, define also $x^{-}_{\mu+\beta\hbar}[\eps]
$ as $\sum_{i\ge 0}\pa^{i}x^{-}_{\mu}[\eps] /\pa\mu^{i}
(\beta\hbar)^{i}/i!$. 

Let us denote by $U^{(e)}_{\la,-}(\tau)$, 
$U^{(f)}_{\la,-}(\tau)$ the subspaces
of $U_{\hbar}\N_{\pm}(\tau)$ linearly
spanned by the products 
$e^{-}_{-\la+2n\hbar}[\eta_{0}] \cdots
e^{-}_{-\la}[\eta_{n}]$, 
resp. by the products 
$f^{-}_{\la}[\eta_{0}] 
\ldots f^{-}_{\la+2n\hbar}[\eta_{n}] $, $n\ge 0$, $\eta_{i}\in k$.   

Let $(\eps_{i})_{i\in \ZZ}$ be a basis of $k$ such that $\eps_{i} = e^{i}$ 
for $i\ge 0$. We can now formulate a Poincar\'e-Birkhoff-Witt result for
$U_{\hbar}\N_{\pm}(\tau)$.

\begin{prop} \label{PBW}
1) Bases of $U^{(e)}_{+}(\tau)$ and of $U^{(f)}_{+}(\tau)$ 
are respectively given by the monomials $(e[\eps_{i_{1}}]\ldots
e[\eps_{i_{p}}])_{0\le i_{1} \le \cdots \le i_{p}}$ and  
$(f[\eps_{i_{1}}]\ldots
f[\eps_{i_{p}}])_{0\le i_{1} \le \cdots \le i_{p}}$; bases of  
$U^{(e)}_{\la,-}(\tau)$ and 
$U^{(f)}_{\la,-}(\tau)$ are respectively given by the 
$(e^{-}_{-\la+2n\hbar}[\eps_{i_{0}}] \cdots
e^{-}_{-\la}[\eps_{i_{n}}])_{i_{0} \le \cdots \le i_{n} < 0}$, 
resp. by the
$(f^{-}_{\la}[\eps_{i_{0}}] 
\ldots f^{-}_{\la+2n\hbar}[\eps_{i_{n}}])_{i_{0} \le \cdots \le i_{n}< 0}$. 

2) The maps 
$U^{(e)}_{\la,-}(\tau) \otimes U^{(e)}_{+} \to U_{\hbar}\N_{+}(\tau)$, 
$U^{(e)}_{+} \otimes U^{(e)}_{\la,-}(\tau) \to U_{\hbar}\N_{+}(\tau)$, 
$U^{(f)}_{\la,-}(\tau) \otimes U^{(f)}_{+} \to U_{\hbar}\N_{-}(\tau)$, 
$U^{(f)}_{+} \otimes U^{(f)}_{\la,-}(\tau) \to U_{\hbar}\N_{-}(\tau)$, 
induced by the multiplication, define vector spaces 
isomorphisms. 
\end{prop}

{\em Proof. \/} We should first derive identities expressing 
the $e[\eps_{i}]e[\eps_{j}], i>j$ 
in terms of combinations of the $e[\eps_{k}]
e[\eps_{l}], k\le l$, for $i,j,k,l\ge 0$;  the 
$e_{-\la+2\hbar}[\eps_{i}]e_{-\la}[\eps_{j}], i>j$ in terms of the 
$e_{-\la+2\hbar}[\eps_{k}]e_{-\la}[\eps_{l}], k\le l$, for $i,j,k,l < 0$; 
the $e_{-\la}[\eps_{i}]e[\eps_{j}],i<0\le j$ 
in terms of the $e[\eps_{k}]e_{-\la}[\eps_{l}]
,l<0\le k$,  
and the $e[\eps_{i}]e_{-\la}[\eps_{j}],j<0\le i$ in terms of the   
 $e_{-\la}[\eps_{k}]e[\eps_{j}],k<0\le l$. 
For this, we may first assume that $\eps_{i} = z^{i}$. Then 
we multiply (\ref{eeps:eeps'}) by $z-w$
and combine the Fourier coefficients as in \cite{Enr-Rub}, sect. 4. 
This proves that the families of 1) generate $U^{(e)}_{+}$ and 
$U^{(e)}_{\la,-}$, and that the two first maps of 2) are surjective. 

The facts that these families are free, and that these maps are 
injective, follow from \cite{Enr-Rub}, Lemma 4.4.
\hfill\qed\medskip

\begin{remark} An informal way to derive (\ref{eeps:eeps'}) is the 
following one. For example, if $\eps=\eps'=+$, we have 
\begin{equation} \label{def-e}
e^{+}_{\la}(z) 
=
\oint_{C_{0}^{-}} 
{{\theta(\zeta - z - \la)}\over{\theta(\zeta - z) \theta(-\la)}} 
e(\zeta) d\zeta , 
\quad
e^{-}_{\la}(z) 
=
\oint_{C_{0,z}} 
{{\theta(\zeta - z - \la)}\over{\theta(\zeta - z) \theta(-\la)}} 
e(\zeta) d\zeta , 
\end{equation}
$C_{0,z}$, resp. $C_{0}^{-}$ being a contour encircling $0$ and $z$ 
(resp. $0$) counterclockwise (resp. clockwise), and 
\begin{equation} \label{def-f}
f^{+}_{ - \la}(z) 
=
\oint_{C_{0}^{-}} 
{{\theta(\zeta - z + \la)}\over{\theta(\zeta - z) \theta(\la)}} 
f(\zeta) d\zeta , 
\quad
f^{-}_{ - \la}(z) 
=
\oint_{C_{0 , z}} 
{{\theta(\zeta - z + \la)}\over{\theta(\zeta - z) \theta(\la)}} 
f(\zeta) d\zeta . 
\end{equation}

Multiply the identity 
$$
{ {\theta(\zeta - \zeta' - \hbar)} \over {\theta(\zeta - \zeta')} } 
e(\zeta) e(\zeta')
=
{ {\theta(\zeta - \zeta' + \hbar)} \over {\theta(\zeta - \zeta')} } 
e(\zeta') e(\zeta)
$$
by ${{\theta(\zeta -z - \la)} \over {\theta(\zeta - z)\theta(-\la)}} 
{{\theta( \zeta' - w - \la)} \over {\theta(\zeta' - w)\theta(-\la)}}$, 
and integrate it over the cycles $C_{0}^{-}$ for $\zeta$, and 
$C_{0}^{\prime -}$ for $\zeta'$, where $C_{0}^{\prime -}$ is a deformation of 
$C_{0}^{-}$, such that $| \zeta | < | \zeta' |$. 
In the resulting identity, replace in the l.h.s., 
$e(\zeta)$ by 
$e^{+}_{\la+\hbar}(\zeta)+ e^{-}_{\la+\hbar}(\zeta)$ and in the r.h.s. by
$e^{+}_{\la-\hbar}(\zeta)+ e^{-}_{\la-\hbar}(\zeta')$. 
The contributions to the integral of the terms in $e^{-}$ 
vanish, 
because these terms are regular at $0$. We then obtain 
\begin{align*}
\oint_{C_{0}^{\prime -}} &   
{ {\theta(\zeta' - w - \la)} \over{\theta(\zeta' - w)\theta( -\la)}}
\lbrace  
{{\theta(z - \zeta' - \hbar)}\over{\theta(z - \zeta')}}
e^{+}_{\la + \hbar}(z)
e^{+}_{\la - \hbar}(\zeta')
\\ & +
{{\theta(\zeta' - z - \la)\theta(-\hbar)}
\over{\theta(\zeta' - z)\theta( - \la)}}
e^{+}_{\la + \hbar}(\zeta')
e(\zeta') \rbrace
d \zeta'
\\ &  
=
\oint_{C_{0}^{\prime -}}
{ {\theta(\zeta' - w - \la)} \over{\theta(\zeta' - w)\theta( -\la)}}
\lbrace
{{\theta(z - \zeta' + \hbar)}\over{\theta(z - \zeta')}}
e^{+}_{\la + \hbar}( \zeta')
e^{+}_{\la - \hbar}(z)
\\ & +
{{\theta(\zeta' - z - \la)\theta(\hbar)}
\over{\theta(\zeta' - z)\theta( - \la)}}
e(\zeta') 
e^{+}_{\la - \hbar}(\zeta')
\rbrace
d \zeta' , 
\end{align*}
that is 
\begin{align} \label{identity}
e^{+}_{\la + \hbar} & (z)  
\lbrace 
{ {\theta(z-w-\la)\theta(\hbar)} \over {\theta(z-w)\theta(-\la)} }
e^{+}_{\la - \hbar}(z)
+
{ {\theta(z-w-\hbar)} \over {\theta(z-w)} }
e^{+}_{\la - \hbar}(w)
\rbrace 
\\ & \nonumber
- 
\lbrace 
{ {\theta(z-w-\la)\theta( - \hbar)} \over {\theta(z-w)\theta(-\la)} }
e^{+}_{\la + \hbar}(z)
+
{ {\theta(z-w+\hbar)} \over {\theta(z-w)} }
e^{+}_{\la + \hbar}(w)
\rbrace 
e^{+}_{\la - \hbar}(z) 
\\ & \nonumber = -
\oint_{C_{0}^{\prime -}}
\theta(-\hbar)
{ {\theta(\zeta'-z-\la)} \over {\theta(\zeta'-z)\theta(-\la)} }
{ {\theta(\zeta'-w-\la)} \over {\theta(\zeta'-w)\theta(-\la)} }
[e^{+}_{\la + \hbar}(\zeta') e(\zeta')
+
e(\zeta') e^{+}_{\la - \hbar}(\zeta') 
]
d \zeta'.  
\end{align}
Specializing (\ref{vertex-e}) for $\zeta = \zeta'$, we find 
$e(\zeta')^{2} = 0$. Therefore 
$$
e^{+}_{\la + \hbar}(\zeta') e(\zeta') +
e(\zeta') e^{+}_{\la - \hbar}(\zeta') = e^{+}_{\la + \hbar}(\zeta')
e^{+}_{\la - \hbar}(\zeta') - e^{-}_{\la + \hbar}(\zeta')
e^{-}_{\la - \hbar}(\zeta').
$$ 
The contribution of the terms in 
$e^{-}$ to the the r.h.s. of (\ref{identity}) is zero, since these 
terms are regular at $0$, 
and the contribution of the terms in $e^{+}$ 
is evaluated by the residues formula.

(Note that this method apparently cannot be adapted 
to derive relations between 
$e^{+}_{\la}(z)$ and $e^{+}_{\la'}(w)$ when $\la \neq \la'\pm 2\hbar$, 
because then the term in $e^{+}$ in the r.h.s. of (\ref{identity})
would no longer be a meromorphic function on $E$.)

As we will see from Thm. \ref{thm:RLL}, 
it is also possible to derive relations between fields
$k^{\pm}(z)$ and $e_{\la}^{\pm}(w), f_{\la}^{\pm}(w)$; for example, we
have 
\begin{equation} \label{k+e+}
k^{+}(z) e_{\la}^{+}(w) k^{+}(z)^{-1}
=
{{\theta(z-w+\hbar)} \over{\theta(z-w)}} e^{+}_{\la+\hbar} (w)
-
{{\theta(z-w-\la)\theta(\hbar)}\over{\theta(z-w)\theta(-\la)}}
e^{+}_{\la + \hbar} (z),
\end{equation}
etc. 
\hfill  \qed
\end{remark} 

\subsection{Shifts in $h$.} \label{ee:sh}

In what follows, we will simply denote $h[1]$ by $h$. 
We will also use the following notation. 
Let us define in each tensor power $A^{\otimes n}$, 
$a^{(i)}$ as the element $1^{\otimes (i-1)} \otimes
a \otimes 1^{\otimes (n-i)}$ for $a\in A$, and for $\beta\in\CC$, 
$x^{-(j)}_{\mu + \hbar\beta  h^{(i)}}[\eps]$ as 
$\sum_{\al\ge 0} {{\pa^{\al}x_{\mu}^{-(j)}[\eps] }\over{\pa \mu^{\al}}}
{{(\hbar \beta h^{(i)})^{\al}}\over{\al!}}$. 

If $n=1$, let us denote also $x^{-(1)}_{\mu + \hbar\beta  h^{(1)}}[\eps]$
simply by $x^{-}_{\mu + \hbar\beta  h}[\eps]$.  Let us also set for
$\eps = +,-$, 
$$
x^{\eps}_{\mu+\hbar \beta h}(z) = \sum_{\al\ge 0}
(\pa/\pa\mu)^{\al}\left( x^{\eps}_{\mu}(z) \right)
(\hbar\beta h)^{\al} / \al! ; 
$$
we have then 
$$
x(z) = x^{+}_{\mu + \hbar \beta h}(z)+ x^{-}_{\mu + \hbar \beta h}(z), 
\quad
x^{-}_{\mu + \hbar\beta  h}(z) =  \sum_{i\in\ZZ}
x^{-}_{\mu + \hbar\beta  h}[\eps^{i}]\eps_{i}(z) .
$$

\begin{lemma} We have the relations
\begin{align} \label{feps:feps':shifts}
&
{ {\theta(z-w+\hbar)} \over {\theta(z-w)} } f^{\eps}_{\la+\hbar h}(z)  
f^{\eps'}_{\la+\hbar h}(w)
+
\eps\eps'{ {\theta(w-z-\la-\hbar h-3\hbar)\theta(\hbar)} 
\over {\theta(w-z)\theta(-\la-\hbar h-3\hbar)} } 
f^{\eps'}_{\la+\hbar h}(w) f^{\eps'}_{\la+\hbar h}(w)
\\ \nonumber & =
{ {\theta(z-w-\hbar)} \over {\theta(z-w)} } f^{\eps'}_{\la+\hbar h}(w) 
f^{\eps}_{\la+\hbar h}(z)
+
\eps\eps'
{ {\theta(z-w-\la-\hbar h -3\hbar)\theta(\hbar)} 
\over {\theta(z-w)\theta(-\la-\hbar h-3\hbar)} } 
f^{\eps}_{\la+\hbar h}(z) f^{\eps}_{\la+\hbar h}(z) ,  
\end{align}
where $\eps,\eps'$ take the values $+$ and $-$. 
\end{lemma}

{\em Proof. \/}
The identity
\begin{align*}
& \sum_{n\ge 0}\phi(\la) (\pa/\pa\la)^{n}
\left( f^{\eps}_{\la}(z_{1}) f^{\eps'}_{\la}(z_{2}) \right)
{{(\hbar h)^{n}}\over{n!}}
\\ & = 
\sum_{p,q,r\ge 0}\phi(\la) {{(\hbar (h+4))^{r}}\over{r!}}
f_{\la}^{\eps}(z_{1}) {{(\hbar (h+2))^{q}}\over{q!}} f_{\la}^{\eps'}(z_{2}) 
{{(\hbar h)^{p}}\over{p!}}
\end{align*}
implies that 
$$
\sum_{n\ge 0}\phi(\la) f_{\la}^{\eps}(z_{1}) f_{\la}^{\eps'}(z_{2}) 
{{(\hbar h)^{n}}\over{n!}}
= 
\phi(\la + \hbar(h+4)) f^{\eps}_{\la+\hbar h + \hbar}(z_{1}) 
f_{\la + \hbar h+ \hbar}^{\eps'}(z_{2}) . 
$$
(\ref{feps:feps':shifts}) then follows from (\ref{feps:feps'}), 
after the change of $\la$ into $\la-\hbar$. 
\hfill \qed

In what follows we will use the notation 
$$
\gamma = -\hbar; 
$$
as we will see, $\gamma$ is connected to the $\eta$ of elliptic quantum
groups by the relation $\gamma = 2 \eta$. 

\subsection{Properties of $K^{+}(z)$.}

Since 
$$
h^{+}(z) = {{\theta'}\over{\theta}}(z) h^{+}[1] + \sum_{i >  0}
h^{+}[e_{i}]e^{i}(z), 
$$
where the $e^{i}$ are elliptic functions, we have 
$$
\left( 
{{q^{\pa} - q^{-\pa}}\over{2\pa}}h^{+} \right)(z) = 
{1\over 2}\ln {{\theta (z + \hbar)}\over{\theta(z - \hbar)}} h^{+}[1]
+\sum_{i > 0} h^{+}[e_{i}] e'_{i}(z), 
$$
with $e'_{i}$ again elliptic functions. 
Therefore, 
\begin{equation} \label{berserkr}
K^{+}(z) = 
\left( {{\theta(z-\gamma )}\over{\theta(z+\gamma)}} \right)^{h/2} 
\sum_{i\ge 0}\wp^{(i)}(z)
\al_{i}
, \quad \al_{i}\in U_{\hbar}\HH(\tau); 
\end{equation}
here $\left( {{\theta(z-\gamma )}\over{\theta(z+\gamma)}} \right)^{h/2}$
is defined as $\exp \left({h\over 2} \ln{{\theta(z+\gamma)}\over
{\theta(z-\gamma)}}\right)$, where the argument of the exponential is
considered as a formal power series in $\gamma$, and we define $\wp$ by 
$\wp = - \pa^{2}(\ln\theta)$. 

\begin{remark}
(\ref{berserkr}) implies that $K^{+}(z)$ has the properties 
\begin{equation} \label{K-per}
K^{+}(z + 1) = K^{+}(z) , \quad
K^{+}(z + \tau) = e^{ - 2 i \pi \hbar h} K^{+}(z);
\end{equation}
informally, we can write 
$$
K^{+}(z) = K^{+}_{ - \gamma h}(z). 
$$
\hfill \qed
\end{remark}

\subsection{Properties of the coproducts} \label{copdt}

\begin{lemma} \label{boker}
Let us fix $\la \in \CC - L$. For $\eps\in k$, we have
\begin{equation} \label{haraldr}
\Delta(e^{-}_{-\la}[\eps]) 
= 
\sum_{i\ge 0} e^{-(1)}_{-\la + \gamma h^{(2)}} [\eps_{i}]
(1\otimes a^{i}_{\la}(\eps)) + e^{-(2)}_{-\la} [\eps] , 
\end{equation}
and 
\begin{equation} \label{visa}
(\Delta \otimes 1) (e^{-(1)}_{-\la + \gamma h^{(2)}}[\eps]) 
= 
\sum_{i\in\ZZ} e^{-(1)}_{-\la + \gamma h^{(2)}
+ \gamma h^{(3)}}[\eps_{i}] \sum_{\al\ge 0}( 1 \otimes 
{{\pa^{\al}{a^{i}_{\la}(\eps)}} \over {\pa \la^{\al}}} \otimes 
{{(\gamma h)^{\al}}\over{\al!}})  
+ 
e^{-(2)}_{-\la + \gamma h^{(3)}}[\eps], 
\end{equation}
where $a^{i}_{\la}$ are linear maps from $k$ to the subalgebra  
$U_{\hbar}\HH_{+}(\tau)$ of $A(\tau)$, generated by the $h^{+}[r],r\in \cO$,
depending holomorphically on $\la\in\CC - L$.  
\end{lemma}

{\em Proof. \/}
We have $\Delta(e(z)) = e(z) \otimes K^{+}(z) + 1\otimes e(z)$, so that 
$$
\Delta(e^{-}_{-\la}[\eps]) = \langle e(z) \otimes K^{+}(z) , p^{-}_{\la}(\eps)
\rangle_{k} 
+ 1\otimes e^{-}_{-\la}[\eps].
$$
By (\ref{berserkr}), we have
$$
\langle e(z) \otimes K^{+}(z) , p^{-}_{\la}(\eps)
\rangle_{k} =
\langle e(z) \otimes 
\left( {{\theta(z-\gamma )}\over{\theta(z+\gamma)}} \right)^{h/2} 
\sum_{i\ge 0}\wp^{(i)}(z) \al_{i} , p^{-}_{\la}(\eps)
\rangle_{k}.  
$$

Note now that the map associating to $\eps\in k$, the series
$$
A_{\la}(\eps) = 
\left( {{\theta(z-\gamma )}\over{\theta(z+\gamma)}} \right)^{h/2} 
\sum_{i\ge 0} {{\pa^{i}p_{\la}^{-}(\eps)} \over{\pa\la^{i}}}
{{(\gamma h)^{i}}\over{i!}}, 
$$
is a linear map from $k$ to $L_{\la}[[\gamma h]]$. 
Therefore 
$$
\left( {{\theta(z-\gamma )}\over{\theta(z+\gamma)}} \right)^{h/2} 
p_{\la}^{-}(\eps) = \sum_{i\ge 0}{{(-\gamma h)^{i}}\over{i!}}
{{\pa^{i}A_{\la}(\eps) }\over{\pa\la^{i}}} , 
$$
and 
$$
\langle e(z) \otimes K^{+}(z) , p^{-}_{\la}(\eps)
\rangle_{k} =
\sum_{i\ge 0} e^{(1)}_{-\la + \gamma h^{(2)}}[\wp^{(i)}A_{\la}(\eps)] 
(1 \otimes \al_{i}). 
$$
(\ref{haraldr}) follows, if we set $a^{i}_{\la}(\eps) = 
\sum_{j\ge 0}\langle \wp^{(j)} A_{\la}(\eps), \rho^{i}\rangle \al_{j}$, 
where $(\rho^{j})_{j\ge 0}$ is the dual basis to $(\eps_{j})_{j\ge 0}$.

(\ref{visa}) then follows directly from (\ref{haraldr}). 
\hfill \qed\medskip

\begin{lemma} \label{erev}
There exists a family of linear maps $(b_{i})_{i\ge 0}$ from $\cO$ to 
the subalgebra $U_{\hbar}\HH(\tau)$ of $A(\tau)$ generated by the 
$h[\eps],\eps\in k$, and $K$, such that 
\begin{equation} \label{Delta:f:+}
\Delta(f[r]) = f[r] \otimes 1 + \sum_{i\ge 0} b_{i}(r)\otimes f[\eps_{i}] , 
\end{equation}
for $r\in \cO$; recall that $(\eps_{i})_{i\ge 0}$ is a basis of $\cO$. 
\end{lemma}

{\em Proof.}
We have $\Delta(f(z)) = f(z) \otimes 1 + q^{- h^{-}(z)} \otimes f(z + \hbar 
K_{1})$, 
so that 
$$
\Delta(f[r]) = f[r] \otimes 1 + \langle  q^{- h^{-}(z)} \otimes f(z + \hbar 
K_{1})
,r\rangle_{k}; 
$$
since $q^{-h^{-}(z)} = \sum_{i\ge 0} \beta_{i}\eps_{i}(z)$, 
for certain $\beta_{i}\in U_{\hbar}\HH(\tau)$, we have
$$
\langle  q^{- h^{-}(z)} \otimes f(z + \hbar K_{1})
,r\rangle_{k} = \langle 
\sum_{i\ge 0} \beta_{i}\eps_{i}(z) \otimes f(z + \hbar K_{1}) , r \rangle_{k}
=
\sum_{i\ge 0} \beta_{i}\otimes f[q^{K_{1}\pa}(r\eps_{i})] . 
$$
We then set for $i\ge 0$, $b_{i}(r) = \sum_{j\ge 0} \beta_{j} \langle
q^{K\pa}(r\eps_{i}) , \rho^{i} \rangle$, where $(\rho^{i})_{i\ge 0}$
is the dual basis to $(\eps_{i})_{i\ge 0}$. 
\hfill \qed \medskip

\begin{lemma} \label{form:bar:Delta}
There exist families of linear maps $(c_{i}^{\la})_{i\ge 0}$ from $\cO$ to 
$U_{\hbar}\HH(\tau)$  and $(d_{i})_{i\in\ZZ}$ from 
$k$ to $U_{\hbar}\HH_{+}(\tau)$, such that the
dependence of $c^{\la}_{i}$ in $\la\in\CC - L$ is holomorphic, and  
\begin{equation} \label{bar:Delta:e+}
\bar\Delta(e[r]) = 1\otimes e[r] + \sum_{i\ge 0}
e[\eps_{i}] \otimes c_{i}^{\la}(r), 
\end{equation}
for $r\in \cO$, 
\begin{equation} \label{bar:Delta:f-}
\bar\Delta(f^{-}_{\la}[\eps]) = f^{-}_{\la}[\eps] \otimes 1
+ \sum_{i\in\ZZ} (d_{i}(\eps) \otimes 1)f^{-(2)}_{\la+\gamma h^{(1)}}
[\eps_{i}]
\end{equation}
and 
\begin{equation} \label{bar:Delta:f-:shift}
\bar\Delta(f^{-(1)}_{\la - \gamma h^{(1)} + 2 \gamma}[\eps]) 
= f^{-(1)}_{\la - \gamma(h^{(1)} + h^{(2)}) + 2 \gamma}[\eps]
+ \sum_{i\in\ZZ} (d_{i}(\eps) \otimes 1)f^{-(2)}_{\la-\gamma h^{(2)} 
+ 2\gamma}
[\eps_{i}].
\end{equation}
\end{lemma}

{\em Proof.\/} 
(\ref{bar:Delta:e+}) is proved in the same way as 
(\ref{Delta:f:+}). Let us prove (\ref{bar:Delta:f-}). We have for 
$\eps\in k$,  
$$
\bar\Delta(f^{-}_{\la}[\eps]) =  f^{-}_{\la}[\eps] \otimes 1 + 
\langle K^{+}(z) \otimes f(z) , p^{-}_{-\la}(\eps)\rangle_{k} ; 
$$
recall that $K^{+}(z) = 
\left( {{\theta(z-\gamma)}\over{\theta(z+\gamma)}}\right)^{h/2} 
\sum_{i\ge 0}
\al_{i} \wp^{(i)}(z)$, with $\al_{i} \in U_{\hbar}\HH_{+}(\tau)$, 
therefore 
$$
\langle K^{+}(z) \otimes f(z) , p^{-}_{-\la}(\eps)\rangle_{k} 
=
\sum_{i\ge 0} \al_{i} \otimes \langle f(z) , \wp^{(i)}(z) p^{-}_{\la}(z)
\left( {{\theta(z-\gamma)}\over{\theta(z+\gamma)}}\right)^{h/2} 
\rangle_{k}; 
$$
since each $\wp^{(i)}(z) p^{-}_{-\la}(z)
\left( {{\theta(z-\gamma)}\over{\theta(z+\gamma)}}\right)^{h/2}$ 
can be expressed as
an expansion 
$$
\sum_{j\ge 0} {{\pa^{j} B_{\la}(z)}\over{\pa\la^{j}}} 
{{(-\gamma h^{(1)})^{j}}\over{j!}},
$$  
with $B_{\la}(z)\in 
L_{-\la}[[\gamma h^{(1)}]]$,
we have 
$$
1\otimes \langle f(z) , \wp^{(i)}(z) p^{-}_{-\la}(z)
{{\theta(z-\gamma h^{(1)})}\over{\theta(z)}} \rangle_{k}
=
f^{-(2)}_{\la + \gamma h^{(1)}} (\la_{i}(\eps)) , 
$$
where $\la_{i}$ are certain linear endomorphisms of $k$. 
This shows (\ref{bar:Delta:f-}). 

(\ref{bar:Delta:f-:shift}) can be deduced from (\ref{bar:Delta:f-}) 
by using the expansion 
$$
f^{-(1)}_{\la - \gamma h^{(1)}}[\eps]
= \sum_{j\ge 0} {{\pa^{j} f^{-(1)}_{\la}[\eps]
}\over{\pa\la^{j}}} {{(-\gamma h^{(1)})^{j}}\over{j!}} , 
$$ 
the identity $\bar\Delta(h) = h^{(1)} + h^{(2)}$ and by 
replacing $\la$ by  $\la + 2 \gamma$. 
\hfill \qed \medskip

\section{Duality} \label{dualty}

Let $U_{\hbar}\G_{+}(\tau)$ be the subalgebra
of $A(\tau)$ generated by $D$, the 
$h[r],r\in \cO$, and $U_{\hbar}\N_{+}(\tau)$, 
and let 
$U_{\hbar}\G_{-}(\tau)$ be the subalgebra
of $A(\tau)$ generated by $K$, the 
$h[\la],\la\in L_{0}$, and 
$U_{\hbar}\N_{-}(\tau)$. 

$(U_{\hbar}\G_{\pm}(\tau),\Delta)$ are Hopf 
subalgebras of $(A(\tau),\Delta)$; 
$(U_{\hbar}\G_{+}(\tau),\Delta)$ and
$(U_{\hbar}\G_{-}(\tau),\Delta')$ are dual
to each other; the duality 
$\langle , \rangle$ is expressed 
by the rules
$$
\langle e[\eps], f[\eps']\rangle
= {1\over \hbar}\langle \eps,\eps' \rangle_{k}, \quad
\langle h[r], h[\la]\rangle
= {2\over \hbar}\langle r,\la \rangle_{k}, \quad
\langle D, K\rangle
= {1\over \hbar}, 
$$ 
the other pairings between generators
being trivial. We denote by 
$\langle, \rangle_{U_{\hbar}\N_{\pm}(\tau)}$
the restriction of $\langle , \rangle$
to $U_{\hbar}\N_{+}(\tau)\times
U_{\hbar}\N_{-}(\tau)$. 

On the other hand, let $U_{\hbar}\bar\G_{+}(\tau)$ and
$U_{\hbar}\bar\G_{-}(\tau)$ be the subalgebras of $A(\tau)$
respectively generated by $U_{\hbar}\N_{+}(\tau)$, $K$ and the $h[\la],
\la\in L_{0}$, and by $U_{\hbar}\N_{-}(\tau)$, $D$ and the $h[r],
r\in \cO$. 
$(U_{\hbar}\bar\G_{\pm}(\tau),\bar\Delta)$ are Hopf 
subalgebras of $(A(\tau),\bar\Delta)$; 
$(U_{\hbar}\bar\G_{+}(\tau),\bar\Delta')$ and
$(U_{\hbar}\bar\G_{-}(\tau),\bar\Delta)$ are dual
to each other; the duality 
$\langle , \rangle'$ is expressed 
by the rules
$$
\langle e[\eps], f[\eps']\rangle'
= {1\over \hbar}\langle \eps,\eps' \rangle_{k}, \quad
\langle h[\la] , h[r] 
\rangle'
= {2\over \hbar}\langle r,\la \rangle_{k}, \quad
\langle K, D\rangle'
= {1\over \hbar}, 
$$ 
the other pairings between generators
being trivial. The restriction of $\langle , \rangle'$ 
to $U_{\hbar}\N_{+}(\tau)$
$\times U_{\hbar}\N_{-}(\tau)$ coincides with
$\langle , \rangle_{U_{\hbar}\N_{\pm}(\tau)}$. 

Let us also denote by $U_{\hbar}\N_{\pm}(\tau)^{[n]}$ the
homogeneous components of degree $n$ (in the $e[\eps]$ or $f[\eps]$) of
$U_{\hbar}\N_{\pm}(\tau)$, and by $U_{\la,-}^{(f);n}$ the intersections 
$U_{\la,-}^{(f)} \cap U_{\hbar}\N_{\pm}(\tau)^{[n]}$. 

Then 

\begin{lemma} \label{orths}
1) The annihilator of $U^{(f)}_{+}$ for 
$\langle , \rangle_{U_{\hbar}\N_{\pm}(\tau)}$ is 
$\sum_{r\in \cO} e[r] U_{\hbar}\N_{+}(\tau)$. 

2) The annihilator of $U^{(e);n}_{\la,-}$ for 
$\langle , \rangle_{U_{\hbar}\N_{\pm}(\tau)}$ is 
$\sum_{\eps\in k} f^{-}_{\la + 2(n-1)\gamma}[\eps]
U_{\hbar}\N_{-}(\tau)^{[n-1]}$.  

3) The annihilator of $U^{(e)}_{+}$ for 
$\langle , \rangle_{U_{\hbar}\N_{\pm}(\tau)}$ is 
$\sum_{r\in \cO} U_{\hbar}\N_{-}(\tau) f[r]$. 

4) The annihilator of $U^{(f);n}_{\la,-}$ for 
$\langle , \rangle_{U_{\hbar}\N_{\pm}(\tau)}$ is 
$\sum_{\eps\in k} U_{\hbar}\N_{+}(\tau)^{[n-1]} 
e^{-}_{-\la +2(n-1)\gamma}[\eps]$. 
\end{lemma}

{\em Proof. \/}
1) and 3) are consequences of \cite{Enr-Rub}, Prop. 6.2. 

Let us show 2). Let us first prove that 
$\sum_{\eps\in k} f^{-}_{\la + 2n\gamma}[\eps]
U_{\hbar}\N_{-}(\tau)^{[n]}$ is
orthogonal to $U^{(e);n+1}_{\la,-}$. Let 
$$
a = e^{-}_{-\la-2n\gamma}[\eta_{n}] \cdots
e^{-}_{-\la}[\eta_{0}], \quad \eta_{i}\in k
$$
belong to $U^{(e);n+1}_{\la,-}$; let $b$ belong to
$U_{\hbar}\N_{-}(\tau)^{[n]}$, $\eps$ belong to $k$ and let us compute
$\langle a,f^{-}_{-\la}[\eps]b \rangle$. This is equal to 
\begin{equation} \label{motti}
\sum_{i}\langle
a_{i}, f^{-}_{-\la}[\eps]\rangle \langle a'_{i}, b  \rangle,
\end{equation} 
where $\Delta(a) = \sum_{i}a_{i} \otimes a'_{i}$. 

{}From (\ref{haraldr}) follows that $\Delta(a)$ is the product of the
terms 
\begin{equation} \label{factor}
\sum_{i\in \ZZ} e^{-(1)}_{-\la -2p\gamma + \gamma h^{(2)}}[\eps_{i}]
(1\otimes a^{i}_{\la + 2 p \gamma})(\eta_{p}) + e^{-(2)}_{-\la - 2 p
\gamma}[\eta_{p}], 
\end{equation}
for $p=n,\ldots,0$. 
This product belongs to $U_{\hbar}\N_{+}(\tau) \otimes
U_{\hbar}\G_{+}(\tau)$. To evaluate (\ref{motti}), we may as well
project the first factor of $\Delta(a)$ on $U_{\hbar}\N_{+}(\tau)^{[1]}$
parallel to all other homogeneous components. 
The contribution of the $(n-p)$th term (\ref{factor}) is 
\begin{align*}
& \sum_{i\in\ZZ}
e^{-(2)}_{-\la - 2 n\gamma}[\eta_{n}]
\cdots e^{-(2)}_{-\la - 2 (p+1)\gamma}[\eta_{p+1}]
e^{-(1)}_{-\la - 2 p\gamma + \gamma h^{(2)}}[\eps_{i}]
(1\otimes a^{i}_{\la + 2p\gamma}(\eta_{p})) 
\\ & 
e^{-(2)}_{-\la - 2 (p-1)\gamma}[\eta_{p-1}]
\cdots
e^{-(2)}_{-\la}[\eta_{0}], 
\end{align*}
that is, using the fact that 
$$
(h-2p) 
e^{-(2)}_{-\la - 2 n\gamma}[\eta_{n}]
\cdots e^{-(2)}_{-\la - 2 (p+1)\gamma}[\eta_{p+1}]
=
e^{-(2)}_{-\la - 2 n\gamma}[\eta_{n}]
\cdots e^{-(2)}_{-\la - 2 (p+1)\gamma}[\eta_{p+1}]
(h-2n), 
$$
 
\begin{align} \label{newterm}
& \sum_{i\in\ZZ}
\sum_{\al\ge 0}
(-\pa/\pa\la)^{\al}(e^{-}_{-\la + 2 n \gamma}[\eps_{i}]) 
\\ & \nonumber \otimes 
{{h^{\al}}\over{\al!}}
e^{-(2)}_{-\la - 2 n\gamma}[\eta_{n}]
\cdots e^{-(2)}_{-\la - 2 (p+1)\gamma}[\eta_{p+1}]
a^{i}_{\la + 2p\gamma}(\eta_{p}) 
e^{-(2)}_{-\la - 2 (p-1)\gamma}[\eta_{p-1}]
\cdots
e^{-(2)}_{-\la}[\eta_{0}] .  
\end{align}
Note now that for any $x\in U_{\hbar}\G_{+}(\tau)$ and $y\in
U_{\hbar}\N_{-}(\tau)$, we have 
\begin{equation} \label{chabat}
\langle hx, y\rangle = 0. 
\end{equation}
Indeed, $\langle hx , y\rangle  = \langle h \otimes x ,
\Delta'(y)\rangle_{(2)}$ 
(denoting by $\langle , \rangle_{(2)}$the tensor square of $\langle ,
\rangle$); but $\Delta'(y)$ belongs to $U_{\hbar} \N_{-}(\tau)
\otimes U_{\hbar}\G_{-}(\tau)$, and $\langle h , U_{\hbar}\N_{-}(\tau)
\rangle = 0$, so that (\ref{chabat}) holds. 

Now the pairing of (\ref{newterm}) with $f^{-}_{\la}[\eps] \otimes b$
is equal to zero either (for $\al = 0$) because $\langle
e^{-}_{-\la}[\eta], f^{-}_{\la}[\eps]\rangle = 0$ for any $\eps,\eta\in
k$ ($L_{\la}$ and $L_{-\la}$ being orthogonal to each other) or by
(\ref{chabat}) for $\al>0$. 

Then standard deformation arguments (see \cite{Enr-Rub}, proof of
Prop. 6.2) show that the orthogonal of 
 $U^{(e);n+1}_{\la,-}$ is exactly 
$\sum_{\eps\in k} f^{-}_{\la + 2 n \gamma}[\eps]
U_{\hbar}\N_{-}(\tau)^{[n]}$.  

Let us now prove 4). Let us first show that for $\eps\in k$,
$U_{\hbar}\N_{+}(\tau)^{[n]}e^{-}_{-\la + 2 n \gamma}[\eps]$ is
orthogonal to $U_{\la,-}^{(f);n+1}$. Let 
$$
a = f^{-}_{\la}[\eta_{0}] \cdots f^{-}_{\la - 2 n \gamma}[\eta_{n}],
\quad \eta_{i}\in k, 
$$
belong to $U_{\hbar}\N_{+}(\tau)^{[n]}$, let $b$ belong to 
$U_{\hbar}\N_{+}(\tau)^{[n]}$, and let us compute $\langle b
e^{-}_{-\la + 2 n \gamma}[\eps], a\rangle$. This is equal to 
\begin{equation} \label{dandolo}
\langle b \otimes e^{-}_{-\la + 2 n \gamma}[\eps],
\bar\Delta(a)\rangle_{(2)}. 
\end{equation} 

{}From (\ref{bar:Delta:f-}) follows that $\bar\Delta(a)$ is the product of
the terms 
\begin{equation} \label{arrigo}
f^{-}_{\la - 2 p \gamma} [\eta_{p}] \otimes 1 + \sum_{i\in\ZZ}
(d_{i}(\eta_{p}) \otimes 1) f^{-(2)}_{\la + \gamma h^{(1)} - 2 p
\gamma}[\eps_{i}],  
\end{equation}
$p = 0,\ldots,n$. Assign degrees $-1$ to terms of the form $f[\eps],
\eps\in k$, and zero to those belonging to $U_{\hbar}\HH(\tau)$; then in
the expansion of the product of the terms (\ref{arrigo}), only those of
degree $-1$ will contribute to (\ref{dandolo}). Therefore (\ref{dandolo})
is equal to 
\begin{align} \label{corner}
& 
\langle b\otimes e^{-}_{-\la + 2 n \gamma}[\eps] , 
\sum_{p=0}^{n} \prod_{k=0}^{p-1}
(f^{-}_{\la - 2 k \gamma}[\eta_{k}]\otimes 1)
\sum_{i\in\ZZ} (d_{i}(\eta_{p}) \otimes 1)
f^{-(2)}_{\la + \gamma ( h^{(1)} - 2 p)}[\eps_{i}]
\\ & \nonumber 
\prod_{k=p+1}^{n}
(f^{-}_{\la - 2 k \gamma}[\eta_{k}]\otimes 1) \rangle_{(2)}. 
\end{align}
Using the identity 
$$
(h-2p) \prod_{k=p+1}^{n}
f^{-}_{\la - 2 k \gamma}[\eta_{k}]
= \left( \prod_{k=p+1}^{n}
f^{-}_{\la - 2 k \gamma}[\eta_{k}] \right) (h-2n) , 
$$
we rewrite (\ref{corner})
as 
\begin{align} \label{murano} 
& 
\langle b\otimes e^{-}_{-\la + 2 n \gamma}[\eps] , 
\sum_{\al\ge 0}
\sum_{p=0}^{n} \sum_{i\in\ZZ}
\left( 
\prod_{k=0}^{p-1} f^{-}_{\la - 2 k \gamma}[\eta_{k}] 
\right)
d_{i}(\eta_{p}) 
\left( 
\prod_{k=p+1}^{n} f^{-}_{\la - 2 k \gamma}[\eta_{k}] 
\right)
{{(\gamma h)^{\al}}\over{\al!}}
\\ & \nonumber 
\otimes (\pa / \pa\la)^{\al} f^{-}_{\la - 2 n \gamma}[\eps_{i}]
\rangle_{(2)} . 
\end{align}
Note now that for $b\in U_{\hbar}\N_{+}(\tau)$, $c\in
U_{\hbar}\bar\G_{-}(\tau)$, 
\begin{equation} \label{kuneitra}
\langle  b, ch \rangle = 0. 
\end{equation}
Indeed, $\langle b , ch \rangle = \langle \bar\Delta'(b),c\otimes h 
\rangle_{(2)} = 0$ because the second factors of the expansion of
$\bar\Delta'(b)$ belong to $U_{\hbar}\N_{+}(\tau)$. Therefore
(\ref{murano}) vanishes, either by (\ref{kuneitra}) or because $\langle
e^{-}_{-\la + 2 n \gamma}[\eps], f^{-}_{\la - 2 n \gamma}[\eps_{i}]
\rangle = 0$.  

It follows that $\sum_{\eps\in k}
U_{\hbar}\N_{+}(\tau)^{[n]}e^{-}_{-\la + 2 n \gamma}[\eps]$ is
orthogonal to $U_{\la,-}^{(f);n+1}$. The same deformation arguments as
above show that these spaces are in fact the orthogonals of each other.
\hfill \qed 

We will also use the following lemma: 

\begin{lemma} \label{F-pair-ell}
1) For $x\in U_{\hbar}\N_{+}(\tau)$, we have 
$$
\langle F, id \otimes x \rangle = x;  
$$
2) for $y\in U_{\hbar}\N_{-}(\tau)$, we have 
$$
\langle F, y \otimes id \rangle = y;  
$$
3) Let $\pi$ be the linear map from $U_{\hbar}\G_{+}(\tau)$ to 
$U_{\hbar}\N_{+}(\tau)$, defined by $\pi(tx) = \varepsilon(t)x$, 
for $x\in U_{\hbar}\N_{+}(\tau)$, $t\in U_{\hbar}\HH_{+}(\tau)$; then we 
have for $x'\in U_{\hbar}\G_{+}$, 
$$
\langle F , id \otimes x' \rangle = \pi(x); 
$$
4) Let $\pi'$ be the linear map from $U_{\hbar}\bar\G_{+}(\tau)$ to 
$U_{\hbar}\N_{+}(\tau)$, defined by $\pi'(yt) = y\varepsilon(t)$, 
for $y\in U_{\hbar}\N_{+}(\tau)$, $t\in U_{\hbar}\HH_{-}(\tau)$; then we 
have for $y'\in U_{\hbar}\G_{+}$, 
$$
\langle F , id \otimes y' \rangle = \pi'(y').  
$$
\end{lemma}

{\em Proof.\/} 1) and 2) are direct consequences of 
\cite{Enr-Rub}, (66) and (68). The proof of 3) and 4) is similar to 
that of Lemma 2.5 of \cite{Rat}. \hfill \qed \medskip

\section{Algebras $A^{+-}$ and $A^{-+}$ and their properties} \label{A+-}

For $X$ a vector space, we denote by $\Hol(\CC - L, X)$ the space 
of holomorphic functions from $\CC - L$ to $X$ and set $\Hol(\CC- L) =
\Hol(\CC - L,\CC)$. 

\begin{defin}
Let us define $A^{-+}$ to be 
the subalgebra of $\Hol(\CC - L , A\otimes A)$
generated (over $\Hol(\CC - L)$) by $h^{(2)}$ and the 
$e^{-(1)}_{-\la + \gamma h^{(2)}}[\eps] f^{(2)}[r]$, with $\eps\in k$
and $r\in \cO$, and $A^{+-}$ as the subalgebra of $\Hol(\CC - L,
A\otimes A)$ 
generated (over $\Hol(\CC - L)$) 
by $h^{(2)}$ and the $e^{(1)}[r]f^{-(2)}_{\la - \gamma h^{(2)} 
+ 2 \gamma}[\eps]$, with $r\in\cO$, $\eps\in k$. 
\end{defin}

\begin{prop} \label{grimr}
The intersection of $A^{+-}$ and $A^{-+}$ is equal to $\Hol(\CC - L, 1\otimes 
\CC[h][[\gamma]])$. 
\end{prop}

{\em Proof. \/} Since we have $[h,f[r]] = -2 f[r]$ for $r\in \cO$, 
we have the relations 
$$
f^{(2)}[r] e^{-(1)}_{-\la + \gamma h^{(2)}}[\eps]
= e^{-(1)}_{-\la + \gamma h^{(2)} +2\gamma}[\eps]f^{(2)}[r],
$$ 
for $\eps\in k,r\in\cO$; therefore 
$A^{-+}$ is linearly spanned by $1$ and the
\begin{equation} \label{xi}
\xi = e^{-(1)}_{-\la + \gamma h^{(2)}}[\eta_{0}]\cdots
e^{-(1)}_{-\la + \gamma h^{(2)}  + 2 n \gamma}[\eta_{n}] 
f^{(2)}[r_{0}]\ldots f^{(2)}[r_{n}] (h^{(2)})^{p}, 
\end{equation}
with $n,p\ge 0$, $\eta_{i}\in k$, $r_{i}\in \cO$. 

On the other hand, 
$A^{+-}$ is linearly spanned by $1$ and the
\begin{equation} \label{eta}
\eta = e^{(1)}[r_{0}] \ldots e^{(1)}[r_{n}] 
f^{-(2)}_{\la - \gamma h^{(2)} + 2\gamma}[\eta_{0}]\cdots
f^{-(2)}_{\la - \gamma h^{(2)} + 2 \gamma}[\eta_{n}] (h^{(2)})^{p},
\end{equation}
$n,p \ge 0$, $\eta_{i}\in k$, $r_{i}\in \cO$. 

Suppose 
that some combination of elements of the form (\ref{eta}) belongs to 
$A^{+-}$. The image of this combination by $l \otimes 1$, 
$l$ any linear form on $A$, is   
some combination 
\begin{equation} \label{comb}
\sum_{p\ge 0}\la_{0;p}(h^{(2)})^{p} 
+ \sum_{n\ge 0}\sum_{i}\la_{i;n;p} 
f^{-}_{\la - \gamma h + 2 \gamma}[\eta_{0}^{(i)}] 
\cdots
f^{-}_{\la - \gamma h + 2 \gamma}[\eta_{n}^{(i)}] (h^{(2)})^{p}, 
\end{equation}
$\la_{0;p},\la_{i;n;p}\in\Hol(\CC - L)$, 
that should belong to $A^{+}$. 

By (\ref{feps:feps'}) with $\eps=\eps'=-$, 
a basis of the linear span of all elements of the form
$$
f^{-}_{\la - \gamma h + 2 \gamma}[\eta_{0}^{(i)}] 
\cdots
f^{-}_{\la - \gamma h + 2 \gamma}[\eta_{n}^{(i)}] (h^{(2)})^{p},
$$
with $\eta_{i}\in k$, is
$$
f^{-}_{\la - \gamma h + 2 \gamma}[\eps_{i_{0}}] 
\cdots
f^{-}_{\la - \gamma h + 2 \gamma}[\eps_{i_{n}}](h^{(2)})^{p},\quad
i_{0}\le \cdots \le i_{n} < 0. 
$$

On the other hand, from (\ref{feps:feps':shifts})
follows that a basis 
of the linear span of all elements of the form 
$f[\eta_{0}]\cdots f[\eta_{n}]h^{p}$, with 
$\eps_{i}\in k$, is given by 
$$
f^{-}_{\la - \gamma h + 2 \gamma}[\eps_{i_0}] \ldots
f^{-}_{\la - \gamma h + 2 \gamma}[\eps_{i_k}]
f[\eps_{i_{k+1}}]
\ldots
f[\eps_{i_{n}}] h^{p} , 
$$
$i_{0} \leq \ldots \leq i_{k} < 0 \leq i_{k+1} \leq \ldots \leq i_{n}$. 

A basis of the intersection of $A^{+}$ with this linear span is
$$
f[\eps_{{i_0}}]
\ldots
f[\eps_{{i_n}}] h^{p}, 
\quad i_{0} \ge \ldots \ge i_{n}\ge 0. 
$$

Therefore the only possibility that (\ref{comb}) belongs to $A^{+}$ is 
that $\la_{i;n;p}=0$ for all $i,n,p$. 
\hfill \qed 
\medskip

\begin{defin}
Let us define $A^{-,\cdot,\cdot}$ as the subspace of the algebra
$\Hol(\CC - L, A^{\otimes 3})$, linearly 
spanned (over $\Hol(\CC - L)$) by the elements of the form 
\begin{equation} \label{medved}
\xi' = e^{-(1)}_{-\la + \gamma ( h^{(2)} + h^{(3)}) 
}[\eta_{1}]\cdots
e^{-(1)}_{-\la + \gamma ( h^{(2)} + h^{(3)}) + 2 (n-1) \gamma
}[\eta_{n}] 
(1\otimes a \otimes b), 
\end{equation}
$n\ge 0$ (recall that the empty product is equal to $1$), 
where $\eta_{i}\in k$, and $a,b\in A$ are 
such that $[h^{(1)} + h^{(2)} + h^{(3)}, \xi'] =0$; 
and $A^{\cdot,\cdot,+}$ as the subspace of $\Hol(\CC - L, A^{\otimes 3})$ 
spanned (over $\Hol(\CC - L)$) by the elements of the form 
$$
\eta' = ( a^{\prime}\otimes b' \otimes 1) 
f^{(3)}[r_{1}]\cdots
f^{(3)}[r_{n}] (h^{(3)})^{s}, \quad n,s\ge 0, 
$$
where $a',b'\in A$, $r_{i}\in \cO$, 
and such that $[h^{(1)} + h^{(2)} + h^{(3)}, \eta'] =0$.
\end{defin}

\begin{prop} \label{egill}
$A^{-,\cdot,\cdot}$ and $A^{\cdot,\cdot,+}$ are subalgebras
of $\Hol(\CC - L, A^{\otimes 3})$. We have
\begin{equation} \label{genet}
(\Delta\otimes 1)(A^{-+}) \subset A^{-,\cdot,\cdot}
\cap A^{\cdot,\cdot,+} , 
\quad
(1\otimes \Delta)(A^{-+}) \subset A^{-,\cdot,\cdot} 
\cap A^{\cdot,\cdot,+}.  
\end{equation}
\end{prop}

{\em Proof. \/} 
That $A^{\cdot,\cdot,+}$ is a subalgebra of $\Hol(\CC - L, A^{\otimes
3})$ follows
easily from its definition. 
Let us show now that $A^{-,\cdot,\cdot}$ is a subalgebra of
$\Hol(\CC - L, A^{\otimes 3})$. Let $\xi'$ and $\xi''$ be elements of
$\Hol(\CC - L, A^{\otimes 3})$ of the 
form (\ref{xi'}), that is 
$$
\xi' = e^{-(1)}_{-\la + \gamma ( h^{(2)} + h^{(3)})}[\eta_{1}]
\cdots
e^{-(1)}_{-\la + \gamma ( h^{(2)} + h^{(3)}) + 2 (n-1) \gamma}[\eta_{n}] 
(1\otimes a \otimes b), 
$$
and
$$
\xi'' = e^{-(1)}_{-\la + \gamma ( h^{(2)} + h^{(3)})}[\eta'_{1}]
\cdots
e^{-(1)}_{-\la + \gamma ( h^{(2)} + h^{(3)}) + 2 (n'-1) \gamma}[\eta'_{n'}] 
(1\otimes a' \otimes b'), 
$$
with $\eta_{i},\eta'_{i} \in k$, $a,b,a',b'\in A$ are such that 
$\xi'$ and $\xi''$ commute with $h^{(1)} + h^{(2)} + h^{(3)}$. 

Since $[h^{(1)}, \xi'] = 2 n \xi$, we have $[h^{(2)} + h^{(3)}, \xi'] 
= - 2 n \xi'$, so that $[h^{(2)} + h^{(3)}, 1 \otimes a \otimes b] = 
- 2n (1\otimes a \otimes b)$. 
It follows that for any $p$, 
$$
(1\otimes a \otimes b)
e^{-(1)}_{-\la + \gamma ( h^{(2)} + h^{(3)}) + 2 (p-1) \gamma}[\eta'_{p}] 
=
e^{-(1)}_{-\la + \gamma ( h^{(2)} + h^{(3)}) + 2 n \gamma + 2(p-1)\gamma}
[\eta'_{p}] 
(1\otimes a \otimes b) .  
$$
The product $\xi'\xi''$ can then be written as 
\begin{align*}
& e^{-(1)}_{-\la + \gamma ( h^{(2)} + h^{(3)})}[\eta_{1}]
\cdots
e^{-(1)}_{-\la + \gamma ( h^{(2)} + h^{(3)}) + 2 (n-1) \gamma}[\eta_{n}] 
\\ & 
e^{-(1)}_{-\la + \gamma ( h^{(2)} + h^{(3)}) + 2 n \gamma}[\eta'_{1}]
\cdots
e^{-(1)}_{-\la + \gamma ( h^{(2)} + h^{(3)}) + 2 (n+n'-1) \gamma}[\eta'_{n'}] 
(1\otimes aa' \otimes bb'), 
\end{align*}
which is of the form (\ref{medved}). Since we also have $[h^{(1)} + 
h^{(2)} + h^{(3)}, \xi'\xi'']=0$, $\xi'\xi''$ belongs to $A^{-,\cdot,\cdot}$.

Let us now prove the first part of (\ref{genet}). 
From (\ref{visa}) follows that for $\eps\in k, r\in \cO$, 
$(\Delta \otimes 1)(e^{-(1)}_{-\la + \gamma h^{(2)}}[\eps] f^{(2)}[r])$ 
is equal to 
$$
\sum_{i\ge 0} \sum_{\al\ge 0}
e^{-(1)}_{-\la + \gamma(h^{(2)}+ h^{(3)})}[\eps_{i}]
{{\pa^{\al} a^{i}_{\la}(\eps)^{(2)}}\over{\pa\la^{\al}}}
{{(\gamma h^{(3)})^{\al}}\over{\al!}} f^{(3)}[r] 
+ e^{-(2)}_{-\la + \gamma h^{(3)}} [\eps] f^{(3)}[r], 
$$
and so belongs to $A^{-,\cdot,\cdot}\cap A^{\cdot,\cdot,+}$. 
We also have $(\Delta \otimes 1)(h^{(2)}) \in A^{-,\cdot,\cdot} \cap 
A^{\cdot,\cdot,+}$.  
Since $h^{(2)}$ and 
the $e^{-(1)}_{-\la - \gamma h^{(2)}} f^{(2)}[r]$, $\eps\in k, 
r\in \cO$, generate $A^{-+}$, and that 
$A^{-,\cdot,\cdot}\cap A^{\cdot,\cdot,+}$ is an algebra, 
we have $(\Delta\otimes 1)(A^{-+}) \subset 
A^{-,\cdot,\cdot}\cap A^{\cdot,\cdot,+}$. 

Let us now prove the second part of (\ref{genet}). Clearly, 
$(1\otimes\Delta)(h^{(2)})$ belongs to $A^{-,\cdot,\cdot}\cap 
A^{\cdot,\cdot,+}$.  
From Lemma \ref{erev} follows that 
for any $\eps\in k,r\in\cO$, $(1\otimes \Delta)(e^{-(1)}_{-\la + 
\gamma h^{(2)}}[\eps]$ $f^{(2)}[r])$ 
is equal to 
$$
e^{-(1)}_{-\la + 
\gamma( h^{(2)} + h^{(3)})}[\eps]\left(f^{(2)}[r]
+ \sum_{i\ge 0}  b_{i}(r)^{(2)} f^{(3)}[\eps_{i}] \right)  , 
$$
and therefore belongs to $A^{-,\cdot,\cdot}\cap A^{\cdot,\cdot,+}$. 
Since $h^{(2)}$ and the 
$e^{-(1)}_{-\la + \gamma h^{(2)}}[\eps] f^{(2)}[r]$, $\eps\in k, 
r\in \cO$, generate $A^{-+}$, and that 
$A^{-,\cdot,\cdot}\cap A^{\cdot,\cdot,+}$ is an algebra, this shows that 
$(1\otimes \Delta)(A^{-+}) \subset A^{-,\cdot,\cdot}\cap A^{\cdot,\cdot,+}$. 
\hfill \qed\medskip

We now define analogues $A^{+,\cdot,\cdot}$ and $A^{\cdot,\cdot,-}$ of
$A^{-,\cdot,\cdot}$ and $A^{\cdot,\cdot,+}$. 

\begin{defin}
$A^{+,\cdot,\cdot}$ is the subspace of the algebra $\Hol(\CC - L,
A^{\otimes 3})$ linearly spanned (over $\Hol(\CC - L)$)
by the elements of the form
\begin{equation} \label{xi'}
\xi' = e^{(1)}[r_{1}] \ldots e^{(1)}[r_{n}] (1\otimes a \otimes b), 
\quad n\ge 0, 
\end{equation}
where $r_{i}\in \cO$, and $a,b\in A$ are such that $[h^{(1)} + 
h^{(2)} + h^{(3)}, \xi']=0$.  
 
$A^{\cdot,\cdot,-}$ is the subspace of $\Hol(\CC - L, A^{\otimes 3})$ 
linearly spanned (over $\Hol(\CC - L)$) by the elements of the form
\begin{equation} \label{eta'}
\eta' = (a'\otimes b'\otimes 1)
f^{-(3)}_{\la - \gamma h^{(3)} + 2 \gamma}[\eta_{1}] \ldots 
f^{-(3)}_{\la - \gamma h^{(3)} + 2 \gamma}[\eta_{n}] (h^{(3)})^{s} , 
\quad n,s\ge 0, 
\end{equation}
where $\eta_{i}\in k$, and $a',b'\in A$ are such that $[h^{(1)} + 
h^{(2)} + h^{(3)}, \eta']=0$.  
\end{defin}

We have then: 

\begin{prop} \label{snorri}
$A^{+,\cdot,\cdot}$ and $A^{\cdot,\cdot,-}$ are subalgebras of 
$\Hol(\CC - L, A^{\otimes 3})$. We have 
\begin{equation} \label{tawil}
(\bar\Delta \otimes 1)(A^{+-}) \subset A^{+,\cdot,\cdot}
\cap A^{\cdot,\cdot,-}, 
\quad
(1\otimes \bar\Delta)(A^{+-}) \subset A^{+,\cdot,\cdot,} \cap 
A^{\cdot,\cdot,-}.  
\end{equation}
\end{prop}

{\em Proof. \/}
That $A^{+,\cdot,\cdot}$  is a subalgebra
of $\Hol(\CC - L, A^{\otimes 3})$ follows easily from its definition. 

Let us now prove that $A^{\cdot,\cdot,-}$ is an algebra. 
Let us consider elements 
of the form (\ref{eta}), 
$$
\eta = (a\otimes b\otimes 1)
f^{-(3)}_{\la - \gamma h^{(3)} + 2 \gamma }[\eta_{1}] \ldots 
f^{-(3)}_{\la - \gamma h^{(3)} + 2 \gamma}[\eta_{n}] (h^{(3)})^{p}
, 
$$
and 
$$
\eta' = (a'\otimes b'\otimes 1)
f^{-(3)}_{\la - \gamma h^{(3)} + 2 \gamma}[\eta'_{1}] \ldots 
f^{-(3)}_{\la - \gamma h^{(3)} + 2 \gamma}[\eta'_{n'}] (h^{(3)})^{p'}, 
$$
where $\eta_{i},\eta'_{i}\in k$, and $a,b,a',b'\in A$ are such that 
$[h^{(1)} + h^{(2)} + h^{(3)}, \eta]=[h^{(1)} + 
h^{(2)} + h^{(3)}, \eta']= 0$.  

Since each
$f^{-(3)}_{\la - \gamma h^{(3)} + 2 \gamma}[\eta_{p}]$ commutes with 
$a'\otimes b'\otimes 1$, 
$\eta\eta'$ can be written as 
\begin{align*}
& \eta\eta' = 
(aa'\otimes bb'\otimes 1) \\ & 
f^{-(3)}_{\la - \gamma h^{(3)} + 2 \gamma}[\eta_{1}] \ldots 
f^{-(3)}_{\la - \gamma h^{(3)} + 2 \gamma}[\eta_{n}]  
f^{-(3)}_{\la - \gamma h^{(3)} + 2 \gamma}[\eta'_{1}] \ldots 
f^{-(3)}_{\la - \gamma h^{(3)} + 2 \gamma}[\eta'_{n'}] (h^{(3)}-2n)^{p}
(h^{(3)})^{p'}, 
\end{align*}
which is of the form (\ref{eta'}); on the other hand, 
$\eta\eta'$ clearly commutes with $h^{(1)} + h^{(2)} + h^{(3)}$, 
which implies that it belongs to $A^{\cdot,\cdot,-}$. 

The proof of (\ref{tawil}) is
similar to that of Prop. \ref{egill}: 
Lemma \ref{form:bar:Delta} implies that for $\eps\in k,r\in\cO$,  
$ (\bar\Delta \otimes 1) (e^{(1)}[r] 
f^{-(2)}_{\la-\gamma h^{(2)}+2\gamma}[\eps]) $  is equal to 
$$
\left( e^{(2)}[r] + \sum_{i\ge 0} e^{(1)}[\eps_{i}] c^{i(2)}_{\la}(r) 
\right) f^{-(3)}_{\la - \gamma h^{(3)} + 2 \gamma}[\eps]
$$
and therefore belongs to $A^{+,\cdot,\cdot}\cap 
A^{\cdot,\cdot,-}$. 
On the other hand, $(1\otimes\bar\Delta)( e^{(1)}[r] 
f^{-(2)}_{\la-\gamma h^{(2)}+2\gamma}[\eps])$ is equal to 
$$
e^{(1)}[r]\left( f^{-(2)}_{\la - \gamma(h^{(2)}+ h^{(3)})+2\gamma}[\eps]
+
\sum_{i\ge 0} d_{i}^{(2)}(\eps) f^{-(3)}_{\la - \gamma h^{(3)}+ 
2\gamma}[\eps_{i}]\right) , 
$$
and so belongs to $A^{+,\cdot,\cdot}\cap A^{\cdot,\cdot,-}$. 
Finally, $(\bar\Delta\otimes 1)(h^{(2)})$ and 
$(1\otimes\bar\Delta)(h^{(2)})$ also belong to $A^{+,\cdot,\cdot}
\cap A^{\cdot,\cdot,-}$.

Since the 
$e^{(1)}[r] f^{-(2)}_{\la-\gamma h^{(2)}+2\gamma}[\eps]$
generate $A^{+-}$, and 
that $A^{+,\cdot,\cdot} \cap A^{\cdot,\cdot,-}$ 
is an algebra, this proves (\ref{tawil}). 
\hfill \qed \medskip 

\begin{prop}
We have 
$$
A^{+,\cdot,\cdot}\cap A^{-,\cdot,\cdot} = \Hol(\CC - L , 1\otimes 
(A^{\otimes 2})^{\HH}), 
\ 
A^{\cdot,\cdot,+}\cap A^{\cdot,\cdot,-} = \Hol(\CC - L , (A^{\otimes 2})^{\HH}
\otimes \CC[h]),  
$$
where $(A^{\otimes 2})^{\HH}$ are the elements of $A^{\otimes 2}$ 
commuting with $h^{(1)} + h^{(2)}$. 
\end{prop}

{\em Proof. \/} The proof is similar to that of Prop. \ref{grimr}. 
\hfill \qed \medskip

\section{Decomposition of $F$.} \label{sect-decomp}

\subsection{Notation.}
We will use the following general notation. 
For $(X_{\la})_{\la\in\CC - L}$ 
a family of maps from an algebra $\cA$ to $A^{\otimes m}$, 
depending on $\la\in\CC - L$
in a holomorphic way, for $n\ge m$ and any injection $k \mapsto i_{k}$ of 
$\{ 1, \ldots, m\}$ into $\{ 1, \ldots, n\}$; for any complex numbers 
$\al_{s}, s = 1,\ldots, n$, we set 
$$
X^{(i_{1}, \ldots, i_{m})}_{\la + \sum_{i=1}^{n}\gamma \al_{k} h^{(k)}} (a)
= 
\sum_{i\ge 0} {{\pa^{i} X_{\la}^{(i_{1},\ldots,i_{k})}}
\over{\pa\la^{i}}} (a) 
{{(\sum_{i=1}^{n}\gamma \al_{k} h^{(k)})^{i}}\over{i!}}, \quad
a\in\cA, 
$$  
where $X_{\la}^{(i_{1},\ldots, i_{k})}(a)$ denotes the image of $X_{\la}(a)$ 
in 
$A^{\otimes m}$ by the map sending the $k$th factor to the $i_{k}$th one. 
If the $X_{\la}$ are algebra morphisms, and that each $\al_{i_{k}}$ vanishes, 
then $X_{\la + \sum_{i=1}^{n} \gamma \al_{k}h^{(k)}}^{(i_{1},
\ldots, i_{k})}$ is an algebra morphism. 

This notation applies in particular if $\cA = \CC$ 
(then $X_{\la}$ is a family of elements of $A^{\otimes m}$). 

\subsection{Decomposition of $F$.}

Let us denote by $U_{\hbar}\N_{\pm}(\tau)^{[n]}$ the homogeneous part of
of $U_{\hbar}\N_{\pm}(\tau)$ of degree $n$. Let us set 
$$
F = \sum_{n\ge 0} F_{n},  \quad \on{with} \quad F_{n}\in
U_{\hbar}\N_{+}(\tau)^{[n]}\otimes U_{\hbar}\N_{-}(\tau)^{[n]}.
$$

\begin{prop} \label{static:decomp}
There exist families $(F^{i;p}_{\la})_{\la\in\CC - L;p\ge 0}$, $i=1,2$
of elements of $A^{\otimes 2}$, 
where $F^{1;p}_{\la}$ is a linear combination with coefficients in
$\Hol(\CC - L)$ of the  
$$
e^{-}_{-\la -2p\gamma}[\eps_{i_{1}}]\cdots
e^{-}_{-\la -2\gamma}[\eps_{i_{p}}]
\otimes f[\eps_{j_{1}}]\cdots f[\eps_{j_{p}}], \quad j_{\al} \ge 0, 
$$
$F^{2;q}_{\la}$ is a similar combination of the 
$$
e[\eps_{i_{1}}]\cdots e[\eps_{i_{q}}] 
\otimes f^{-}_{\la + 2 q \gamma}[\eps_{j_{1}}]\cdots 
f^{-}_{\la + 2 \gamma}[\eps_{j_{q}}], \quad i_{\al}\ge 0, 
$$
and $F^{i;0}_{\la}=1,i=1,2$, such that 
\begin{equation} \label{Decomp}
F_{n} = \sum_{p+q = n} F^{2;q}_{\la + 2 p \gamma} F^{1;p}_{\la} . 
\end{equation}
\end{prop}

{\em Proof. \/} Let us define the linear maps $\Pi^{(e)}_{+,\la}$ and 
$\Pi^{(e)}_{-,\la}$, from $U_{\hbar}\N_{+}(\tau)$ to $U^{(e)}_{+}$, 
resp. to $U^{(e)}_{-,\la}$, by 
$$
\Pi^{(e)}_{+,\la}(ab) = a\varepsilon(b), \quad 
\Pi^{(e)}_{-,\la}(ab) = \varepsilon(a)b, 
\quad a \in U^{(e)}_{+}, b\in U^{(e)}_{-,\la}, 
$$
and 
$\Pi^{(f)}_{+,\la}$ and 
$\Pi^{(f)}_{-,\la}$, from $U_{\hbar}\N_{-}(\tau)$ to $U^{(f)}_{+}$, 
resp. to $U^{(f)}_{-,\la}$, by 
$$
\Pi^{(f)}_{+,\la}(ab) = \varepsilon(a)b, \quad 
\Pi^{(f)}_{-,\la}(ab) = a\varepsilon(b), 
\quad  a\in U^{(f)}_{-,\la}, b\in U^{(f)}_{+}. 
$$
Note that $\Pi^{(e)}_{+,\la}$ is a left $U^{(e)}_{+}$-module map,
and $\Pi^{(f)}_{+,\la}$ is a right $U^{(f)}_{+}$-module map. 
From Prop. \ref{PBW} also follows that the kernels of  
$\Pi^{(e)}_{+,\la}$, $\Pi^{(e)}_{-,\la}$ 
$\Pi^{(f)}_{+,\la}$ and of $\Pi^{(f)}_{-,\la}$ 
are respectively  
$\sum_{\eps\in k} U_{\hbar}\N_{+}(\tau) e^{-}_{-\la}[\eps]$, 
$\sum_{r\in \cO} e[r]U_{\hbar}\N_{+}(\tau)$, 
$\sum_{\eps\in k} f^{-}_{\la}[\eps]U_{\hbar}\N_{-}(\tau)$ and 
$\sum_{r\in \cO} U_{\hbar}\N_{-}(\tau)f[r]$.

\begin{lemma} 
1) For each $n\ge 0$, 
$(1\otimes \Pi_{+,\la + 2 n \gamma}^{(f)}) (F_{n})$ 
and 
$(\Pi_{-,\la + 2 \gamma}^{(e)} \otimes 1)(F_{n})$ 
both belong to $U^{(e)}_{\la+2\gamma,-} \otimes U^{(f)}_{+}$; 

2) for each $n\ge 0$, 
$(1\otimes \Pi_{-,\la + 2 n \gamma}^{(f)}) (F_{n})$ and 
$(\Pi_{+,\la + 2 \gamma}^{(e)} \otimes 1) (F_{n})$
both belong to $U^{(e)}_{+}\otimes U^{(f)}_{\la+2n\gamma,-}$;

3) we have the equalities 
$$
(1\otimes \Pi_{+,\la + 2 n \gamma}^{(f)}) (F_{n})
= (\Pi_{-,\la + 2 \gamma}^{(e)} \otimes 1)(F_{n})
$$
and
$$
(\Pi_{+,\la + 2 \gamma}^{(e)} \otimes 1) (F_{n})
= (1\otimes \Pi_{-,\la + 2 n \gamma}^{(f)}) (F_{n}). 
$$ 
\end{lemma} 

{\em Proof of the lemma. \/} Let us prove the first part of 1). 
$(1\otimes\Pi^{(f)}_{+,\la+2n\gamma})(F_{n})$ 
clearly belongs to $U_{\hbar}\N_{+}(\tau) \otimes U_{+}^{(f)}$. 
Let now $x$ belong to $\sum_{\eps\in k} f^{-}_{\la + 2 n\gamma}[\eps] 
U_{\hbar}\N_{-}(\tau)^{[n-1]}$. Consider 
$\langle (1\otimes\Pi^{(f)}_{+,\la+2n\gamma})(F_{n}) , x \otimes id \rangle$; 
this is equal to $\Pi_{+,\la + 2n\gamma}(x)$ by Lemma \ref{F-pair-ell}, 1), 
and therefore to zero. 
By Lemma \ref{orths}, 2), it follows that 
$(1\otimes\Pi^{(f)}_{+,\la+2n\gamma})(F_{n})$  also belongs to 
$U_{\la + 2\gamma,-}^{(e)} \otimes U_{\hbar}\N_{-}(\tau)$. 

The proof of the second part of 1) and of 2) is similar, and uses the
other statements of Lemma \ref{orths}, and the above description of the 
kernels of $\Pi^{(e)}_{-,\la+2\gamma}$, $\Pi^{(f)}_{-,\la+2n\gamma}$, 
and $\Pi^{(e)}_{+,\la+2\gamma}$. 

Let us now prove the first part of 3). Let us fix $a_{+}$ in 
$U_{+}^{(f);n}$ and $a_{-}$ in $U^{(e);n}_{\la+2\gamma,-}$. 
Then 
\begin{align*} 
& \langle (1\otimes \Pi^{(f)}_{+,\la+2n\gamma})(F_{n}) - 
(\Pi^{(e)}_{-,\la+2\gamma}\otimes 1)(F_{n}) 
, a_{+} \otimes a_{-} 
\rangle_{(2)}
\\ & 
= \langle a_{-}, \Pi^{(f)}_{+,\la+2n\gamma}(a_{+}) \rangle
- \langle \Pi^{(e)}_{-,\la+2\gamma}(a_{-}), a_{+}\rangle
\\ & = \langle a_{-}, a_{+}\rangle - \langle a_{-}, a_{+}\rangle 
=0.
\end{align*}
Since $(1\otimes \Pi^{(f)}_{+,\la+2n\gamma})(F_{n}) - 
(\Pi^{(e)}_{-,\la+2\gamma}\otimes 1)(F_{n})$ 
belongs to $U^{(e);n}_{\la+2\gamma,-} \otimes U^{(f);n}_{+}$, 
and that the pairing $\langle , \rangle$ induces injections 
from $U^{(e);n}_{\la+2\gamma,-}$ 
in the dual of $U^{(f);n}_{+}$ 
and from $U^{(f);n}_{+}$ 
in the dual of $U^{(e);n}_{\la+2\gamma,-}$,
$(1\otimes \Pi^{(f)}_{+,\la+2n\gamma})(F_{n}) - 
(\Pi^{(e)}_{-,\la+2\gamma}\otimes 1)(F_{n})$ 
is equal to zero.

The proof of the second part of 3) is similar. 
\hfill \qed \medskip

Let us set now 
\begin{equation} \label{F1lan}
F_{\la}^{1;n} = (1\otimes \Pi_{+,\la + 2 n \gamma}^{(f)}) (F_{n})
= (\Pi_{-,\la + 2 \gamma}^{(e)} \otimes 1)(F_{n})
\end{equation}
and
$$
F_{\la}^{2;n} = (\Pi_{+,\la + 2 \gamma}^{(e)} \otimes 1) (F_{n})
= (1\otimes \Pi_{-,\la + 2 n \gamma}^{(f)}) (F_{n}). 
$$ 

To prove (\ref{decomp:ell}), let us consider the family (indexed by $\la
\in \CC - L$) of linear endomorphisms 
$\ell_{\la}$ of $U_{\hbar}\N_{+}(\tau)$, defined by 
$$
\ell_{\la}(x) = \langle \sum_{p+q=n} F^{2;q}_{\la+2p\gamma}F^{1;p}_{\la}, 
id \otimes x\rangle, 
$$ 
for $x$ in $U_{\hbar}\N_{+}(\tau)^{[n]}$. 

Assign in $U_{\hbar}\G_{+}(\tau)$, the degree $0$ to the elements of  
$U_{\hbar}\HH_{+}(\tau)$, and the degree $1$ to each $e[\eps],\eps\in k$. 
Let us denote by $U_{\hbar}\G_{+}(\tau)^{[q]}$ the subspace of 
$U_{\hbar}\G_{+}(\tau)$ formed by the elements of degree $q$. 
Then 
$$
\Delta(U_{\hbar}\N_{+}(\tau)^{[n]}) \subset 
\sum_{p+q=n}U_{\hbar}\N_{+}(\tau)^{[p]}
\otimes U_{\hbar}\G_{+}(\tau)^{[q]}. 
$$

Let us now fix $x$ in $U_{\hbar}\N_{+}(\tau)^{[n]}$, and let us set 
$\Delta(x) = \sum_{i,p+q=n} x'_{p,i} \otimes x''_{q,i}$, with 
$x'_{p,i} \in U_{\hbar}\N_{+}(\tau)^{[p]}$, 
$x''_{q,i} \in U_{\hbar}\G_{+}(\tau)^{[q]}$. 
Expand $\ell_{\la}(x)$ as 
\begin{align} \label{expr-ell1}
\ell_{\la}(x) & \nonumber = 
\sum_{i, p+q = n} 
\langle F^{2;q}_{\la+2p\gamma}, id\otimes x'_{q,i} \rangle 
\langle F^{1;p}_{\la} , id \otimes 
x''_{p,i} \rangle 
\\ & \nonumber = \sum_{i, p+q = n}
\Pi^{(e)}_{+,\la + 2 \gamma + 2p\gamma} (x'_{q,i})
\Pi^{(e)}_{-,\la + 2 \gamma} (\langle F , id \otimes x''_{q,i}\rangle)
\\ & = \sum_{i,p+q=n} 
\Pi^{(e)}_{+,\la + 2 \gamma + 2p\gamma} (x'_{q,i})
(\Pi^{(e)}_{-,\la + 2 \gamma} \circ\pi) (x''_{q,i}); 
\end{align}
the first equality follows from the Hopf algebra pairing rules,
the second one from Lemma \ref{F-pair-ell}, 1), 2), and the third one 
from the same Lemma, 3). 
This formula enables us to show: 

\begin{lemma} \label{cara}
$\ell_{\la}$ is a left $U^{(e)}_{+}(\tau)$-module map. 
\end{lemma}

{\em Proof. \/} Recall that the product map defines a linear isomorphism 
from 
$U_{\hbar}\HH_{+}(\tau)\otimes U^{(e)}_{+} \otimes 
U^{(e)}_{\la+2\gamma,-}$ 
onto $U_{\hbar}\G_{+}(\tau)$. 
Define $U_{\hbar}\B_{+}(\tau)$ as the image of 
$U_{\hbar}\HH_{+}(\tau)\otimes U^{(e)}_{+}\otimes 1$  
by this map. 
$U_{\hbar}\B_{+}(\tau)$ is then a subalgebra of $U_{\hbar}\G_{+}(\tau)$. 

On the other hand, since for $r\in \cO$, $\Delta(e[r])$ is equal to 
$\langle e(z) \otimes K^{+}(z) , r\rangle_{k} + 1 \otimes e[r]$,  
it belongs to 
$U_{\hbar}\N_{+}(\tau) \otimes U_{\hbar}\B_{+}(\tau)$. 
It follows that $\Delta(U_{+}^{(e)}) \subset U_{\hbar}\N_{+}(\tau) \otimes 
U_{\hbar}\B_{+}(\tau)$.

$\Pi^{(e)}_{-,\la+2\gamma}\circ\pi$ is now defined as follows: to $x
\in U_{\hbar}\G_{+}(\tau)$ decomposed as 
$$
\sum_{i} h_{i}x^{+}_{i}x^{-}_{i},
\quad \on{with} \quad  h_{i} \in U_{\hbar}
\HH_{+}(\tau), 
 x^{+}_{i}\in U^{(e)}_{+},
x^{-}_{i}\in U^{(e)}_{\la+2\gamma,-},
$$ 
it associates 
$\Pi^{(e)}_{-,\la+2\gamma}(\sum_{i}\varepsilon(h_{i})x^{+}_{i}x^{-}_{i})$, 
that is $\sum_{i} \varepsilon(h_{i}x^{+}_{i})x^{-}_{i}$. Therefore, it 
satisfies 
\begin{equation} \label{module3}
(\Pi^{(e)}_{-,\la+2\gamma} \circ \pi)(bx) =
\varepsilon(b)(\Pi^{(e)}_{-,\la+2\gamma} \circ \pi)(x) ,  
\end{equation}
for $x\in U_{\hbar}\G_{+}(\tau)$, $b\in U_{\hbar}\B_{+}(\tau)$. 

Let us fix now $x \in U_{\hbar}\N_{+}(\tau)^{[n]}$, 
$b\in U^{(e),m}_{+}$.
Let us set $\Delta(x) = \sum_{i,p+q=n} x'_{p,i} \otimes x''_{q,i}$, 
$x'_{p,i}, x''_{q,i}$ as above, and $\Delta(b) = \sum_{j,p'+q' = m} b'_{p',j} 
\otimes b''_{q',j}$, $b'_{p',j}\in U_{\hbar}\N_{+}(\tau)^{[p']}$, 
$b''_{q',j}\in U_{\hbar}\B_{+}(\tau)^{[q']}$ 
(where $U_{\hbar}\B_{+}(\tau)^{[q']}$ is the intersection of 
$U_{\hbar}\B_{+}(\tau)$ and $U_{\hbar}\G_{+}(\tau)^{[q']}$).

Then 
\begin{align*} 
\ell_{\la}(bx) & = 
\sum_{i,j,p+q=n,p'+q'=m} \Pi^{(e)}_{+,\la+2\gamma+2(p+p')\gamma}
(b'_{q',j}x'_{q,i}) (\Pi^{(e)}_{-,\la+2\gamma}\circ\pi)(b''_{p',j}x''_{p,i})
\\ & = 
\sum_{i,j,p+q=n,p'+q'=m} \Pi^{(e)}_{+,\la+2\gamma+2(p+p')\gamma}
(b'_{q',j}x'_{q,i}) \varepsilon(b''_{p',j})
(\Pi^{(e)}_{-,\la+2\gamma}\circ\pi)(x''_{p,i}) , 
\end{align*}
where the first equality follows from (\ref{expr-ell1}), and the second one
from (\ref{module3}). 
If $p'\neq 0$, $\varepsilon(b''_{p',j})$ vanishes, so that
$\sum_{j} b'_{m,j}\varepsilon(b''_{0,j}) = \sum_{j,p'+q'=m}
b'_{q',j}\varepsilon(b''_{p',j}) = b$. Therefore $\ell_{\la}(bx)$ is equal to 
\begin{align*}
& \sum_{i,j,p+q=n} \Pi^{(e)}_{+,\la+2\gamma+2p\gamma}
(bx'_{q,i}) 
(\Pi^{(e)}_{-,\la+2\gamma}\circ\pi)(x''_{p,i})
\\ &  
 = \sum_{i,j,p+q=n} b\Pi^{(e)}_{+,\la+2\gamma+2p\gamma}
(x'_{q,i}) 
(\Pi^{(e)}_{-,\la+2\gamma}\circ\pi)(x''_{p,i}) = b\ell_{\la}(x), 
\end{align*}
where the first equality follows from the fact that 
$\Pi^{(e)}_{-,\la+2\gamma}$ is a left $U^{(e)}_{+}$-module map. 
\hfill \qed \medskip

Set now for $x$ in $U_{\hbar}\N_{+}(\tau)^{[n]}$, 
$\bar\Delta(x) = \sum_{i,p+q=n} \bar x'_{p,i}\otimes \bar x''_{q,i}$, 
with $\bar x'_{p,i}\in U_{\hbar}\N_{+}(\tau)^{[p]}$, $\bar x''_{q,i}
\in U_{\hbar}\bar \G_{+}(\tau)^{[q]}$. 
Then 
\begin{align} \label{expr-2}
\nonumber 
\ell_{\la}(x) & = \sum_{i,p+q = n}
\langle F^{2;q}_{\la+2p\gamma}, id \otimes \bar x''_{q,i} \rangle'
\langle F^{1;p}_{\la}, id\otimes \bar x'_{p,i} \rangle'
\\ & \nonumber
= \sum_{i,p+q = n}
\Pi^{(e)}_{+,\la+2p\gamma+2\gamma}(\langle F,id\otimes \bar x''_{q,i}\rangle') 
\Pi^{(e)}_{-,\la+2\gamma}(\bar x'_{p,i})
\\ & = \sum_{i,p+q = n}
(\Pi^{(e)}_{+,\la+2p\gamma+2\gamma}\circ\pi') (\bar x''_{q,i}) 
\Pi^{(e)}_{-,\la+2\gamma}(\bar x'_{p,i}), 
\end{align}
using the Hopf algebra pairing rules, and Lemma \ref{F-pair-ell}, 3). 
Before we use this identity to derive a result analogous to Lemma 
2.7 of \cite{Rat}, we will show the following results. 

\begin{lemma} \label{ghagen}
For any $\eps,\eps'\in k$, the product $(q^{-h^{-}})[\eps]
e^{-}_{\la}[\eps']$ is equal to some combination 
$\sum_{i}e^{-}_{\la-2\gamma}[\eta'_{i}](q^{-h^{-}})[\eta_{i}]$, for certain
$\eta_{i},\eta'_{i}\in k$. 
\end{lemma}

{\em Proof.\/} Recall that in (\ref{K-:e}), the ratio of theta-functions 
should be expanded for the argument of $K^{-}$ near $0$; therefore, we
have 
$$
(q^{-h^{-}})(z) e(w) (q^{h^{-}})(z)
=
\sum_{\al\ge 0} (\pa / \pa w)^{\al}\left( {{\theta(w-\gamma K - \gamma)}
\over{\theta(w-\gamma K +\gamma)}}\right) 
{{(-z)^{\al}}\over{\al!}} e(w). 
$$
It follows that for $\eps'\in k$, we have
$$
(q^{-h^{-}})(z) e[\eps'] (q^{h^{-}})(z)
=
\sum_{\al\ge 0} {{(-z)^{\al}}\over{\al!}}
e\left[ \eps'(w)
(\pa / \pa w)^{\al}\left( {{\theta(w-\gamma K - \gamma)}
\over{\theta(w-\gamma K + \gamma)}}\right) \right] , 
$$
therefore we have for $\eps\in k$, 
$$
(q^{-h^{-}}) [\eps] e[\eps'] 
=
\sum_{\al\ge 0} 
e\left[ \eps' (w)
(\pa / \pa w)^{\al}\left( {{\theta(w-\gamma K - \gamma)}
\over{\theta(w-\gamma K + \gamma)}}\right) \right] 
(q^{h^{-}}) [\eps(z) {{(-z)^{\al}}\over{\al!}}] . 
$$

Now if $\eps'$ belongs to $L_{\la}$, the products 
$\eps'(w)
(\pa / \pa w)^{\al}\left( {{\theta(w-\gamma K - \gamma)}
\over{\theta(w-\gamma K + \gamma)}}\right)$
belong to $L_{\la-2\gamma}$. The lemma follows. 
\hfill \qed\medskip 

\begin{lemma} \label{weak-mod}
For $x\in U_{\hbar}\N_{+}(\tau)$, $\eps\in k$, we have: 

1) $\Pi^{(e)}_{-,\la}(x e^{-}_{-\la}[\eps]) = 
\Pi^{(e)}_{\la}(x) e^{-}_{-\la}[\eps]$; 

2) $\Pi^{(e)}_{+,\la}(x e^{-}_{-\la}[\eps]) = 0$. 
\end{lemma}

{\em Proof.\/} This follows directly from the definitions of 
$\Pi^{(e)}_{\pm,\la}$. \hfill \qed\medskip

\begin{lemma} \label{vaggio}
We have the identities 
$$
\ell_{\la-2\gamma}(x e^{-}_{-\la}[\eps]) = 
\ell_{\la}(x) e^{-}_{-\la}[\eps], 
$$
for $x\in U_{\hbar}\N_{+}(\tau), \eps\in k$. 
\end{lemma}

{\em Proof. \/}
Let us fix $x\in U_{\hbar}\N_{+}(\tau)^{[n]}$, and set as above
$\bar\Delta(x) = \sum_{i,p+q=n} \bar x'_{p,i}\otimes \bar x''_{q,i}$, 
with $\bar x'_{p,i}\in U_{\hbar}\N_{+}(\tau)^{[p]}$, $\bar x''_{q,i}
\in U_{\hbar}\bar \G_{+}(\tau)^{[q]}$. 
From formula (\ref{bar:Delta:e:ell}) follows that 
$\bar x''_{q,i}$ can be decomposed as a sum $\sum_{j}y_{q,i,j}h_{q,i,j}$, 
where $y_{q,i,j}$ belongs to $U_{\hbar}\N_{+}(\tau)^{[q]}$ and $h_{q,i,j}$
is a linear combination of products of the form $(q^{-h^{-}})[\eps_{1}]
\cdots (q^{-h^{-}})[\eps_{p}]$, with $\eps_{i}\in k$ (where we denote 
$\langle q^{-h^{-}(z)}, \eta(z)\rangle_{k}$ by $(q^{-h^{-}})[\eta]$). 

Let now $\eps$ belong to $k$; we have 
$$
\bar\Delta'(xe^{-}_{-\la}[\eps]) = \sum_{i,j,p+q = n}
y_{q,i,j} h_{q,i,j} \otimes \bar x'_{p,i}
\left( \langle q^{-h^{-}(z)} \otimes e(z), p_{-\la}(\eps)
\rangle_{k} 
+ e^{-}_{-\la}[\eps]\otimes 1 \right) . 
$$
Let us set 
$$
\langle q^{-h^{-}(z)} \otimes e(z), p_{-\la}(\eps) \rangle_{k}
= \sum_{s\in\ZZ} (q^{-h^{-}})[\eps_{s}] \otimes e[\rho_{s}(\eps)], 
$$
where $\rho_{s}$ is a family of linear endomorphisms of $k$. 

The according to (\ref{expr-2}), we have
\begin{align} \label{yad}
\ell_{\la-2\gamma} (x e^{-}_{-\la}[\eps]) & = 
\sum_{s,i,j,p+q=n} \left( \Pi^{(e)}_{+,l\la+2p\gamma+2\gamma}
\circ\pi'\right) (y_{q,i,j}h_{q,i,j}(q^{-h^{-}})[\eps_{s}]) 
\Pi^{(e)}_{-,\la}(\bar x'_{p,i} e[\rho_{s}(\eps)])
\\ & \nonumber + 
\sum_{s,i,j,p+q=n} \left( \Pi^{(e)}_{+,l\la+2p\gamma}
\circ\pi'\right) (y_{q,i,j}h_{q,i,j} e^{-}_{-\la}[\eps]) 
\Pi^{(e)}_{-,\la}(\bar x'_{p,i}) . 
\end{align}

Using the property of $\pi'$ that $\pi'(xt) = \pi'(x) \varepsilon(t)$, 
for $x\in U_{\hbar}\G_{-}(\tau)$, $t\in U_{\hbar}\HH_{-}(\tau)$, we write 
the first sum of the r.h.s. of (\ref{yad}) as 
\begin{align} \label{nua}
\sum_{s,i,j,p+q=n} & \left( \Pi^{(e)}_{+,\la+2p\gamma+2\gamma}
\circ\pi'\right) (y_{q,i,j}h_{q,i,j}) \varepsilon((q^{-h^{-}})[\eps_{s}]) 
\Pi^{(e)}_{-,\la}(\bar x'_{p,i} e[\rho_{s}(\eps)])
\\ & \nonumber 
= \sum_{i,j,p+q=n} \left( \Pi^{(e)}_{+,\la+2p\gamma+2\gamma}
\circ\pi'\right) (y_{q,i,j}h_{q,i,j}) 
\Pi^{(e)}_{-,\la}(\bar x'_{p,i} e^{-}_{-\la}[\eps])
\\ & \nonumber 
= \sum_{i,j,p+q=n} \left( \Pi^{(e)}_{+,\la+2p\gamma+2\gamma}
\circ\pi'\right) (y_{q,i,j}h_{q,i,j}) 
\Pi^{(e)}_{-,\la}(\bar x'_{p,i}) e^{-}_{-\la}[\eps]
\\ & = \ell_{\la}(x) e^{-}_{-\la}[\eps]; 
\end{align}
here the first equality follows from the properties of $\varepsilon$, 
and the second one from Lemma \ref{weak-mod}, 1). 

According to Lemma \ref{ghagen}, each product $h_{q,i,j}e^{-}_{-\la}[\eps]$ 
can be written as a sum 
$$
\sum_{t\in\ZZ} e^{-}_{-\la-2p\gamma}[\eps_{t}]
h_{q,i,j,t},
$$   
with $h_{q,i,j,t}\in U_{\hbar}\HH_{-}(\tau)$.
It follows that the second sum of the r.h.s. of (\ref{yad}) can be written
as
\begin{equation} \label{gur}
\sum_{t,s,i,j,p+q=n}
\left( \Pi^{(e)}_{+,\la+2p\gamma}\circ\pi'\right) 
(y_{q,i,j}e^{-}_{-\la-2p\gamma}[\eps_{t}]h_{q,i,j,t}) 
\Pi^{(e)}_{-,\la}(\bar x'_{p,i}) . 
\end{equation} 
But 
\begin{align*}
\left( \Pi^{(e)}_{+,\la+2p\gamma} \circ\pi'\right) &  
(y_{q,i,j}e^{-}_{-\la-2p\gamma}[\eps_{t}]h_{q,i,j,t}) 
\\ & =
\Pi^{(e)}_{+,\la+2p\gamma}(y_{q,i,j}e^{-}_{-\la-2p\gamma}[\eps_{t}]) 
\varepsilon(h_{q,i,j,t}) 
 = 0
\end{align*}
by Lemma \ref{weak-mod}, 2). 
Therefore (\ref{gur}) vanishes. The lemma follows from this and 
(\ref{nua}). 
\hfill \qed \medskip

We are now in position to show that for any $\la\in \CC - L$, $x\in U_{\hbar}
\N_{+}$, 
\begin{equation} \label{zac}
\ell_{\la}(x) = x. 
\end{equation}
Using Prop. \ref{PBW}, decompose $x$ as a sum 
$\sum_{i,p,q} x^{+}_{p,i}x^{-}_{q,i}$, with 
$x^{+}_{p,i}$ in $U^{(e);p}_{+}$, $x^{-}_{q,i}$ in 
$U^{(e);q}_{\la,-}$. By Lemma \ref{cara}, $\ell_{\la}(x)$ is equal to 
$\sum_{i,p,q} x^{+}_{p,i} \ell_{\la} (x^{-}_{q,i})$; and by Lemma 
\ref{vaggio}, this last expression is equal to 
$\sum_{i,p,q} x^{+}_{p,i} \ell_{\la+2q\gamma}(1) x^{-}_{q,i}$. 
We easily check that for any $\la\in\CC - L$, $\ell_{\la}(1)=1$. 
(\ref{zac}) follows. 

The proposition now follows from the comparison of 
(\ref{zac}) and Lemma \ref{F-pair-ell}, and from the fact that 
$\langle , \rangle_{U_{\hbar}\N_{\pm}(\tau)}$ is non-degenerate.  
\hfill \qed \medskip

We can now obtain another decomposition of $F$:  

\begin{corollary} \label{dyn:decomp}
There is a unique a decomposition of $F$ as 
\begin{equation} \label{decomp:ell}
F= F^{2}_{\la}F^{1}_{\la} , \quad \on{with} \quad  
F^{1}_{\la}\in A^{-+}\quad \on{and} \quad 
F^{2}_{\la}\in A^{+-},  
\end{equation}
with $(\varepsilon\otimes 1)(F_{i}) = (1\otimes \varepsilon)(F_{i}) = 1$, 
$i=1,2$. 
\end{corollary}

{\em Proof. \/}
Set 
$$
F^{2}_{\la} = 
\sum_{q\ge 0} \sum_{\al\ge 0} 
(\pa / \pa\la)^{\al}(F^{2;q}_{\la}) 
{{(-\gamma h^{(2)})^{\al}}\over{\al!}},
$$
and 
\begin{equation} \label{F1la}
F^{1}_{\la} =  
\sum_{p\ge 0} \sum_{\al\ge 0} 
(\pa / \pa\la)^{\al}(F^{1;p}_{\la}) 
{{(-\gamma h^{(2)})^{\al}}\over{\al!}} ; 
\end{equation}
$F^{1}_{\la}$and $F^{2}_{\la}$ belong repectively to $A^{-+}$ and
$A^{+-}$. Since we have also 
$$
F^{2}_{\la} = 
\sum_{q\ge 0} \sum_{\al\ge 0} 
\pa_{\la}^{\al}(F^{2;q}_{\la + 2 p \gamma}) 
{{(-\gamma(h^{(2)}+2p))^{\al}}\over{\al!}},
$$ 
we can write 
$$
F^{2}_{\la} F^{1}_{\la} = \sum_{p,q\ge 0}
\sum_{\al\ge 0} \pa_{\la}^{\al} (F^{2;q}_{\la + 2 p \gamma}
F^{1;p}_{\la})
{{(-\gamma h^{(2)})^{p}}\over{p!}} = F. 
$$

Let us now prove the unicity of the decomposition (\ref{decomp:ell}). 
Let $(F^{\prime 1}_{\la} , F^{\prime 2}_{\la})$ some other solution to 
(\ref{decomp:ell}). Then, by Prop. \ref{grimr}, we will have
$F^{\prime 2}_{\la} = F^{2}_{\la} u$,   $F^{\prime 1}_{\la} = 
u^{-1}F^{1}_{\la}$, with $u$ some invertible element of 
$\Hol(\CC - L, 1 \otimes \CC[h][[\gamma]])$. On the other hand, 
$(\varepsilon \otimes 1)(F^{2}_{\la}) = (\varepsilon \otimes 1)
(F^{\prime 2}_{\la}) = 1$ implies that $u=1$. 
\hfill \qed \medskip

\begin{lemma} \label{approx}
We have an expansion 
$$
F^{1}_{\la} \in 1 + \hbar \sum_{i\ge 0} e^{-(1)}_{-\la + \gamma h^{(2)}}
[e_{i;-\la}] f^{(2)}[e^{i}]+ U_{\hbar}\N_{+}(\tau)^{\ge 2}
\otimes U_{\hbar}\N_{-}(\tau)^{\ge 2}\CC[h], 
$$
where $U_{\hbar}\N_{\pm}(\tau)^{\ge 2} = \oplus_{i\ge 2}
U_{\hbar}\N_{\pm}(\tau)^{[i]}$. 

\end{lemma}

{\em Proof. \/} This follows from formulas (\ref{F1la}), (\ref{F1lan}), 
and from the fact that $\Pi_{-,\la+2\gamma}^{(e)}$ maps each 
$U_{\hbar}\N_{+}(\tau)^{[i]}$ to itself. 
\hfill \qed

\section{Twisted cocycle property.}

Let us define for $\la\in \CC - L$, 
$$
\Phi_{\la} = F_{\la-\gamma h^{(3)}}^{1(12)}
(\Delta\otimes 1)(F_{\la}^{1})
\left( F_{\la}^{1(23)} (1\otimes \Delta)(F_{\la}^{1})\right)^{-1} . 
$$

\begin{prop} The family 
$(\Phi_{\la})_{\la\in\CC - L}$ belongs to $ A^{-,\cdot,\cdot}
\cap A^{\cdot,\cdot,+}$.
\end{prop}

{\em Proof. \/}
First observe that if $(\phi_{\la})_{\la\in\CC - L}$ belongs to 
$A^{-+}$, then $(\phi_{\la - \gamma h^{(3)}}^{(12)})_{\la\in\CC - L}$ and
$(\phi_{\la}^{(23)})_{\la\in\CC - L}$ belong to 
$A^{-,\cdot,\cdot}\cap A^{\cdot,\cdot,+}$;
this follows easily from the definitions of these algebras. 
It follows that $(F_{\la - \gamma h^{(3)}}^{1(12)})_{\la\in\CC - L}$ and
$(F_{\la}^{1(23)})_{\la\in\CC - L}$ also belong to 
$A^{-,\cdot,\cdot}\cap A^{\cdot,\cdot,+}$.

By (\ref{genet}), the families
$(\Delta\otimes 1)(F^{1}_{{\la}})$ and $(1\otimes
\Delta)(F^{1}_{{\la}})$ also belong to $A^{-,\cdot,\cdot} \cap
A^{\cdot,\cdot,+}$. Since $A^{-,\cdot,\cdot} \cap A^{\cdot,\cdot,+}$ 
is an algebra, and has the following property: any $x\in A^{-,\cdot,\cdot}
\cap A^{\cdot,\cdot,+}$, invertible in $\Hol(\CC - L, A^{\otimes 3})$, 
is such that $x^{-1}$
belongs to $A^{-,\cdot,\cdot}\cap A^{\cdot,\cdot,+}$, 
$(\Phi_{\la})_{\la\in\CC - L}$ 
belongs to $A^{-,\cdot,\cdot}\cap A^{\cdot,\cdot,+}$. 
\hfill \qed\medskip 

Using (\ref{F-cocycle}), we may rewrite $\Phi_{\la}$ as 
$$
\Phi_{\la} = \left((\bar\Delta\otimes 1)(F^{2}_{\la})(F^{2(12)}_{\la -
\gamma h^{(3)}}) \right)^{-1} (1\otimes
\bar\Delta)(F_{\la}^{2}) 
(1\otimes F^{2(23)}_{\la}) . 
$$

\begin{prop}
$\Phi_{\la} \in A^{+,\cdot,\cdot} 
\cap A^{\cdot,\cdot,-}$.
\end{prop}

{\em Proof. \/} We now remark that if $(\psi_{\la})_{\la\in\CC - L}$ 
belongs to $A^{+-}$, 
then $(\psi_{\la -\gamma h^{(3)}}^{(12)})_{\la\in\CC - L}$ 
and $(\psi_{\la}^{(23)})_{\la\in\CC - L}$ 
both belong to $A^{+,\cdot,\cdot}\cap A^{\cdot,\cdot,-}$.
It follows that $(F_{\la -\gamma h^{(3)}}^{2(12)})_{\la\in\CC - L}$ 
and $(F_{\la}^{2(23)})_{\la\in\CC - L}$ 
also belong to $A^{+,\cdot,\cdot}\cap A^{\cdot,\cdot,-}$.

By (\ref{tawil}), the families $(\bar\Delta\otimes 1)(F_{\la}^{2})$ and 
$(1\otimes \bar\Delta)(F_{\la}^{2})$ also belong to 
$A^{+,\cdot,\cdot}\cap A^{\cdot,\cdot,-}$.
Since $A^{+,\cdot,\cdot}\cap A^{\cdot,\cdot,-}$
also has the property
that any $x\in A^{+,\cdot,\cdot}
\cap A^{\cdot,\cdot,-}$, invertible in $\Hol(\CC - L, A^{\otimes 3})$, 
is such that $x^{-1}$
belongs to $A^{+,\cdot,\cdot}\cap A^{\cdot,\cdot,-}$, 
$(\Phi_{\la})_{\la\in\CC - L}$ 
belongs to $A^{+,\cdot,\cdot}\cap A^{\cdot,\cdot,-}$. 
\hfill \qed\medskip 

{}From the two above propositions follows that we have 
\begin{equation} \label{form-of-Phi}
\Phi_{\la}  = \sum_{i\ge 0}1 \otimes a^{(i)}_{\la} \otimes h^{i}, 
\end{equation}
for a certain family $a^{(i)}_{\la}$ of elements of $A(\tau)$, commuting with 
$h$. 

Let us define now 
$$
\Delta_{\la} = \Ad(F^{1}_{\la}) \circ \Delta; 
$$
this is a family of algebra morphisms from $A(\tau)$ to
$A(\tau)^{\otimes 2}$, depending
on $\la\in\CC - L$ in a holomorphic way. 

Then we have the twisted quasi-Hopf condition
\begin{equation} \label{tw-q-H}
(\Delta_{\la - \gamma h^{(3)}}\otimes 1)\circ 
\Delta_{\la} = \Ad(\Phi_{\la})\circ (1\otimes
\Delta_{\la})\circ\Delta_{\la}. 
\end{equation}

\begin{prop} \label{compat} 
$(\Delta_{\la})_{\la\in\CC - L}$ and $(\Phi_{\la})_{\la\in\CC - L}$ 
satisfy the compatibility condition
\begin{align} \label{maoz}
& (\Delta_{\la - \gamma (h^{(3)} + h^{(4)})} \otimes 1 \otimes 1)
(\Phi_{\la})
(1\otimes 1 \otimes \Delta_{\la}) (\Phi_{\la})
\\ & \nonumber
=
\Phi^{(123)}_{\la - \gamma h^{(4)}} (1\otimes \Delta_{\la - \gamma h^{(4)}}
\otimes 1)
(\Phi_{\la}) \Phi_{\la}^{(234)}. 
\end{align}
\end{prop}

{\em Proof. \/} The proof is a straightforward computation. 
We can understand it in the following way. 

For $(V_{i},\rho_{i})_{i=1,2}$ two $A(\tau)$-modules (that is, we are given 
algebra morphisms $\rho_{i}$ from $A(\tau)$ to $\End(V_{i})$), define the
family of $A(\tau)$-modules $(V_{1}\otimes V_{2}, \rho_{1}\otimes_{\la} 
\rho_{2})_{\la\in\CC - L}$ (sometimes simply denoted by $V_{1}\otimes_{\la} 
V_{2}$) by $\rho_{1}\otimes_{\la}\rho_{2} = (\rho_{1}\otimes\rho_{2})
\circ \Delta_{\la}$. 

Then for $(V_{i},\rho_{i})_{i=1,2,3}$ three $A(\tau)$-modules, (\ref{tw-q-H}) 
implies that the image of 
$\Phi_{\la}$ by $\otimes_{i=1}^{3}\rho_{i}$ is an intertwiner of 
$A(\tau)$-modules from $V_{1} \otimes_{\la} (V_{2}\otimes_{\la}V_{3})$ to 
$(V_{1}\otimes_{\la - \gamma h^{(3)}} V_{2})\otimes_{\la} V_{3}$.

For $(V_{i},\rho_{i})_{1\le i \le 4}$ four $A(\tau)$-modules, we then have the 
sequences of $A(\tau)$-modules morphisms 
\begin{align*}
& 
V_{1} \otimes_{\la} (V_{2}\otimes_{\la} (V_{3}\otimes_{\la} V_{4}))
\to
V_{1} \otimes_{\la} ( (V_{2}\otimes_{\la-\gamma h^{(4)}} 
V_{3})\otimes_{\la} V_{4})
\\ & \to
(V_{1} \otimes_{\la-\gamma h^{(4)}} (V_{2}\otimes_{\la-\gamma h^{(4)}} 
V_{3}))\otimes_{\la} V_{4}
\to 
((V_{1} \otimes_{\la-\gamma(h^{(3)}+h^{(4)})} V_{2})
\otimes_{\la-\gamma h^{(4)}} V_{3})\otimes_{\la} V_{4}
\end{align*}
given by the image by $\otimes_{i=1}^{4}\rho_{i}$ of the r.h.s. of 
(\ref{maoz}), and the sequence
\begin{align*}
& V_{1} \otimes_{\la} (V_{2}\otimes_{\la} (V_{3}\otimes_{\la} V_{4}))
\to
(V_{1} \otimes_{\la-\gamma(h^{(3)}+h^{(4)})} V_{2}) \otimes_{\la}
(V_{3} \otimes_{\la} V_{4})
\\ & \to
((V_{1} \otimes_{\la-\gamma(h^{(3)}+h^{(4)})} V_{2})
\otimes_{\la-\gamma h^{(4)}} V_{3})\otimes_{\la} V_{4}
\end{align*}
given by the image by $\otimes_{i=1}^{4}\rho_{i}$ of the l.h.s. of the
same equation; these sequences of morphisms coincide.  
\hfill \qed\medskip

We are ready to conclude: 

\begin{thm} \label{tw-cocycle}
We have $\Phi_{\la} = 1$. Therefore 
\begin{equation} 
F^{1(12)}_{\la - \gamma h^{(3)}} (\Delta\otimes 1)(F^{1}_{\la})
=
F^{1(23)}_{\la} (1\otimes \Delta)(F^{1}_{\la}).
\end{equation}  
\end{thm}

{\em Proof.\/}
After substitution of (\ref{form-of-Phi}), the l.h.s. of (\ref{maoz}) 
becomes the product of 
$$
\sum_{i\ge 0} 1\otimes 1 \otimes a_{\la}^{(i)}\otimes h^{i}
$$
and 
$$
\sum_{i\ge 0} 1\otimes a_{\la}^{(i)} \otimes (h\otimes 1 + 1 \otimes h)^{i}; 
$$
since $h$ commutes with the $a^{(i)}_{\la}$, it follows that these two terms 
commute. Therefore (\ref{maoz}) simplifies to 
$$
(1\otimes 1 \otimes \Delta_{\la})(\Phi_{\la})
=
(\Phi_{\la - \gamma h^{(4)}} \otimes 1) (1\otimes 
\Delta_{\la - \gamma h^{(4)}} \otimes 1) (\Phi_{\la}). 
$$

Apply now $1\otimes \varepsilon \otimes 1 \otimes 1$ to 
this identity. Since $(1\otimes \varepsilon\otimes 1)(\Phi_{\la}) = 1$,
the first and second term map to $1$. On the other hand, since 
$(\varepsilon \otimes 1)\circ \Delta_{\la} = id$, the last term maps
to $\Phi_{\la}$. Therefore, $\Phi_{\la} = 1$. 
\hfill \qed\medskip

\begin{remark} It would be interesting to find some analogue for the
$M(\la)$ of 
\cite{BBB}, that is some family of elements of $A(\tau)$ whose twisted
coboundary would be $F_{\la}$. 
This element should belong to some completion of $A(\tau)$, and
for $\gamma = 0$ coincide with a ``longest Weyl group element the affine 
algebra''.  
\end{remark}

\section{Dynamical Yang-Baxter equation}

The paper \cite{BBB}, sect. 2, contains the following result:  

\begin{prop} (see \cite{BBB}) Let $(\cA,\Delta_{\infty}^{\cA},
\cR_{\infty}^{\cA})$ 
be a quasi-triangular Hopf algebra, with a fixed element $\wt h$. 
Let $F(\la)$ be a family of invertible elements of $\cA \otimes \cA$, 
parametrized by some subset $U \subset \CC$. 
Set $\Delta (\la) = \Ad(F(\la)) \circ 
\Delta^{\cA}_{\infty}$.  Suppose that the identity
\begin{equation} \label{tw-coc}
F^{(12)}(\la-\gamma \wt h^{(3)}) (\Delta_{\infty}^{\cA}\otimes
1)(F(\la))
=
F^{(23)}(\la) (1\otimes \Delta_{\infty}^{\cA})(F(\la))
\end{equation}
is satisfied. Then we have 
\begin{equation} \label{tw:coass}
(\Delta (\la - \gamma \wt h^{(3)}) \otimes 1) \circ \Delta (\la)
=
(1 \otimes \Delta (\la)) \circ \Delta (\la) ,
\end{equation}
and  if we set $\cR(\la) = F^{(21)}(\la)\cR_{\infty}^{\cA}F(\la)^{-1}$, 
we have the identity
\begin{equation} \label{DYBE}
\cR^{(12)}(\la) \cR^{(13)}(\la - \gamma \wt h^{(2)}) \cR^{(23)}(\la) 
=
\cR^{(23)}(\la - \gamma \wt h^{(1)}) \cR^{(13)}(\la) \cR^{(12)}(\la - 
\gamma \wt h^{(3)}) . 
\end{equation}
\end{prop}

{\em Proof. \/} Define $\wt \Delta^{[12]}(\la)$ as the linear map from
$\cA^{\otimes 2}$ to $\cA^{\otimes 3}$ defined by 
\begin{equation} \label{Del12}
\wt \Delta^{[12]}(\la)(x\otimes y) = F^{(12)}(\la) 
(\Delta_{\infty}^{\cA}
\otimes 1)(x\otimes y) F^{(12)}(\la - \gamma \wt h^{(3)})^{-1},  
\end{equation}
and $\wt \Delta^{[21]}(\la)(x\otimes y) = \wt \Delta^{[12]}(\la)(x\otimes
y)^{(213)}$. Then we have 
\begin{equation} \label{QT}
\cR^{(12)}(\la) 
\wt \Delta^{[12]}(\la)(x)
= \wt \Delta^{[21]}(\la)(x) \cR^{(12)}(\la - \gamma \wt h^{(3)}),  \quad
x\in \cA^{\otimes 2}
\end{equation}
and
\begin{equation} \label{Del12R}
\wt \Delta^{[12]}(\la) (\cR(\la)) = \cR^{(13)}(\la-\gamma \wt h^{(2)})
\cR^{(23)}(\la). 
\end{equation}
Applying (\ref{QT}) to $x = \cR(\la)$ yields (\ref{DYBE}). 

We could also define $\wt \Delta^{[23]}(\la)$ by 
$$
\wt \Delta^{[23]}(x\otimes y) = F^{(32)}(\la - \gamma \wt h^{(1)})
(1\otimes \Delta^{\cA \prime}_{\infty})(x\otimes y) F^{(32)}(\la)^{-1}, 
$$
$\wt\Delta^{[32]}(T) = \wt\Delta^{[23]}(T)^{(132)}$, for $T\in
\cA^{\otimes 2}$; then we have  
$\wt \Delta^{[23]}(\la)(x)\cR^{(23)}(\la) =
\cR^{(23)}(\la-\gamma \wt h^{(1)}) \wt \Delta^{[32]}(T)$, and 
$$
\wt \Delta^{[23]}(\la) (\cR(\la)) = \cR^{(12)}(\la)\cR^{(13)}(\la-\gamma
\wt h^{(2)}). 
$$
Note also the identity $\wt \Delta^{[23]}(T) = (\wt
\Delta^{[12]}(T^{(21)-1}))^{(312)-1}$.  
\hfill\qed\medskip 

Identities (\ref{tw-coc}), (\ref{tw:coass}) and (\ref{DYBE})
are respectively called the twisted cocycle condition for the family 
$F(\la)$, the twisted coassociativity condition for $\Delta(\la)$,
and the dynamical Yang-Baxter equation for $\cR(\la)$. 

\begin{thm} \label{Thm:DYBE}
Let us set in $A(\tau)^{\otimes 2}$, 
$$
\cR_{\infty} =  q^{D \otimes K} q^{{1\over 2}\sum_{i\ge 0} h[e^{i}] 
\otimes h[e_{i;0}]} 
q^{\sum_{i\in\ZZ} e[\eps^{i}] \otimes f[\eps_{i}]}
$$
and for $\la\in \CC - L$, $\cR_{\la} = (F^{1}_{\la})^{(21)} \cR_{\infty}
(F^{1}_{\la})^{-1}$. Then the family
$(\cR_{\la})_{\la\in \CC - L}$ satisfies the dynamical Yang-Baxter
relation 
\begin{equation} \label{DYBE:ell}
\cR^{(12)}_{\la} \cR^{(13)}_{\la - \gamma h^{(2)}} \cR^{(23)}_{\la} 
=
\cR^{(23)}_{\la - \gamma h^{(1)}} \cR^{(13)}_{\la} \cR^{(12)}_{\la - 
\gamma h^{(3)}}. 
\end{equation}
\end{thm}

{\em Proof. \/} This follows directly from the above proposition and the
fact that $( A(\tau) , \Delta , $ $\cR_{\infty})$ 
is a quasi-triangular Hopf
algebra (see \cite{Enr-Rub}).  
\hfill \qed \medskip 

\begin{remark}
Let $\wt \Delta_{\la}^{[12]}$ be the linear map from $A^{\otimes 2}$ to
$A^{\otimes 3}$ defined by (\ref{Del12}), and $\Delta^{[12]}_{\la}$ the
map $\wt \Delta^{[12]}_{\la} \otimes id$ from $A^{\otimes 2} \otimes
\Diff(\CC - L)$ to $A^{\otimes 3} \otimes \Diff(\CC - L)$. Here
$\Diff(\CC - L)$ is the ring of differential operators in $\la\in\CC -
L$. Let us set 
$$
\cR = e^{\gamma h^{(1)}\pa_{\la}}\cR_{\la}. 
$$
Then we have the simple relation $\Delta^{[12]}_{\la}(\cR) =
\cR^{13}\cR^{23}$. The form of $\cR$ indicates that the algebra element
$h[{\theta' \over \theta}]$ can be naturally added with the derivative
${\pa \over {\pa\la}}$. This indication could be useful for the study
of twisted conformal blocks: in such a theory we need to add
differential operators to the elements of $\G_{\la}$. 
\end{remark}

\section{Level $0$ representations of $A(\tau)$, $L$-operators and
$RLL$ relations}

In \cite{HGQG}, we studied the $2$-dimensional representations, at level
$0$, of the quantum groups introduced there. In the case of the algebra
$A(\tau)$, these representations can be described as follows. 

Let us denote by $k_{\zeta}$ the local field $\CC((\zeta))$, by
$\pa_{\zeta}$ its derivation $d/d\zeta$, and by $k_{\zeta}[\pa_{\zeta}]$
the associated ring of differential operators. 
Let $(v_{1},v_{-1})$ be the standard 
basis of $\CC^{2}$, and $E_{ij}$ the endomorphism of $\CC^2$ 
defined by $E_{ij}(v_{\alpha}) =\delta_{\alpha, j} v_{i}$. 

\begin{prop}  \label{fd:repres} (see \cite{HGQG}, Prop. 9)
There is a morphism of algebras 
$\pi_{\zeta} : A(\tau) \to \End(\CC^{2}) \otimes k_{\zeta}[\pa_{\zeta}]
[[\gamma]]$, defined by the formulas 
$$
\pi_{\zeta}(K) = 0, \quad \pi_{\zeta}(D) = \Id_{\CC^{2}} \otimes \pa_{\zeta}, 
$$
$$
\pi_{\zeta}(h[r]) = 
E_{11}  \otimes \left( {{2}\over{1+q^{\pa}}} r \right)
(\zeta) 
- E_{-1-1} \otimes \left( {{2} \over {1+q^{-\pa}}} r \right) (\zeta) ,
\quad
r\in \cO, 
$$
$$
\pi_{\zeta}(h[\la]) = E_{11} \otimes \left( {{ 1 - q^{-\pa}}\over{\hbar
\pa }} \la \right) (\zeta) - E_{-1-1} \otimes 
\left( {{ q^{\pa} - 1} \over {\hbar\pa}} \la \right) (\zeta) , 
\quad \la\in L_{0}, 
$$
$$
\pi_{\zeta}(e[\eps]) = { { \theta(\hbar)} \over \hbar } 
E_{1,-1}\otimes \eps(\zeta), 
\quad 
\pi_{\zeta}(f[\eps]) = E_{-1,1}\otimes \eps(\zeta), 
\quad \eps\in k.  
$$
\end{prop}

\begin{lemma} \label{adar}
The image of $\cR_{\la}$ by $\pi_{\zeta} \otimes \pi_{\zeta'}$ is 
\begin{equation} \label{im-R}
(\pi_{\zeta} \otimes \pi_{\zeta'})(\cR_{\la-\gamma}) = 
A(\zeta,\zeta') R^{-}( \zeta-\zeta', \la) , 
\end{equation}
where 
\begin{align} \label{R-mat} 
R^{-}(z,\la) & =  E_{11} \otimes E_{11} + E_{-1,-1} \otimes E_{-1,-1} 
+ 
{{\theta(z)}\over {\theta(z + \gamma)}}  
E_{11} \otimes E_{-1,-1} 
\\ & \nonumber 
+ { {\theta(\la - \gamma) \theta(\la + \gamma)} 
\over{\theta(\la)^{2}}} 
{{\theta(z)}\over{\theta(z+\gamma)}}  
E_{-1,-1} \otimes E_{11} 
+ { {\theta(z+\la) \theta(\gamma)}  
\over {\theta(z+\gamma) \theta(\la)}}
E_{1,-1} \otimes E_{-1,1} 
\\ & \nonumber  
-
{{\theta(z-\la)\theta(\gamma)} 
\over{\theta(z+\gamma)\theta(\la)}}
E_{-1,1} \otimes E_{1,-1} , 
\end{align} 
and 
$A(\zeta,\zeta')$ is equal to $\exp( 
\sum_{i\ge 0}\left({1\over
\pa}{{q^{\pa}-1}\over{q^{\pa}+1}} e^{i}\right)(\zeta)
e_{i;0}(\zeta'))$. 
\end{lemma}

{\em Proof. \/} Since the image by $\pi_{\zeta}$ and $\pi_{\zeta'}$
of $U_{\hbar}\N_{\pm}(\tau)^{\ge 2}$ is zero, and by Lemma
\ref{approx}, this image is the same as that of 
$$
(1 + \hbar \sum_{i\in \ZZ} e^{-(2)}_{-\la+\gamma h^{(1)}}[\eps_{i}]
f^{(1)}[\eps^{i}]) q^{D \otimes K}
q^{{1\over 2} \sum_{i\in\ZZ}h[e^{i}] \otimes h[e_{i;0}]}
(1 + \hbar \sum_{i\ge 0} e^{(1)}[\eps^{i}] f^{-(2)}_{\la - \gamma h^{(2)} +
2 \gamma} [\eps_{i}]) . 
$$
After we use the expansions 
$$
\sum_{i\ge 0} e_{i;\la}(z) e^{i}(w) = {{\theta(z-w+\la)}
\over{\theta(z-w)\theta(\la)}}, \on{\ for\ } \la\in\CC - L, 
\quad
\sum_{i\ge 0} e_{i;0}(z) e^{i}(w) = {{\theta'}\over{\theta}}(z-w), 
$$
and the identities 
$$
\sum_{i\ge 0} (f(\pa)e^{i})(\zeta) e_{i;0}(\zeta') = 
\sum_{i\ge 0} e^{i}(\zeta) (f(-\pa)e_{i;0})(\zeta'), 
$$
for $f$ any polynomial in $\pa$, and
$$
\exp\left( {{q^{\pa_{z}} - 1}\over{\pa_{z}}}
{{\theta'}\over{\theta}}(z-w)\right) =  
{{\theta(z-w+\hbar)}\over{\theta(z-w)}}
, 
$$
we find 
\begin{align*} 
(\pi_{\zeta} \otimes \pi_{\zeta'})(\cR_{\la}) & = 
\\ \nonumber &  
A(\zeta,\zeta')
\left( 1 + \theta(\hbar) (E_{-1,1} \otimes E_{1,-1}) 
{ {\theta(\zeta' - \zeta + \la + \gamma)}
\over{\theta(\zeta' - \zeta) \theta(\la+\gamma)}} \right)
\cdot 
\\ & \nonumber
\left( 
E_{1,1}\otimes E_{1,1} + E_{-1,-1} \otimes E_{-1,-1}
\right. 
\\ & \nonumber \left. + 
{{\theta(\zeta'-\zeta)}\over{\theta(\zeta'-\zeta+\hbar)}}
E_{1,1}\otimes E_{-1,-1} + 
{{\theta(\zeta'-\zeta-\hbar)}\over{\theta(\zeta'-\zeta)}} 
E_{-1,-1} \otimes E_{1,1}
\right) \cdot  
\\ & \nonumber
\left( 1 - \theta(\hbar) (E_{1,-1} \otimes E_{-1,1}) 
{ {\theta(\zeta - \zeta' + \la + \gamma)}
\over{\theta(\zeta - \zeta') \theta(\la + \gamma)}} \right) ;   
\end{align*}
the lemma follows. 
\hfill \qed \medskip 

Define $R^{+}(z,\la)$ as $R^{-}(z,\la)^{-1}$. We then have 
\begin{align} \label{R+}
R^{+}(z,\la) & = E_{11} \otimes E_{11} + E_{-1,-1} \otimes E_{-1,-1}
+ {{\theta(z)}\over{\theta(z-\gamma)}}
{{\theta(\la-\gamma)\theta(\la+\gamma)}\over{\theta(\la)^{2}}}
E_{1,1} \otimes E_{-1,-1} \\ & \nonumber
+ 
{{\theta(z)}\over{\theta(z-\gamma)}} E_{-1,-1} \otimes E_{11} 
- {{\theta(z+\la)\theta(\gamma)}\over {\theta(z-\gamma)\theta(\la)}} 
E_{1,-1} \otimes E_{-1,1}
\\ & \nonumber
+ {{\theta(z-\la)\theta(\gamma)}\over{\theta(z-\gamma)\theta(\la)}} 
E_{-1,1} \otimes E_{1,-1} . 
\end{align}

Let us define now the $L$-operators as follows.  
Set 
$$
L_{\la}^{+}(\zeta) = (1\otimes \pi_{\zeta})(\cR_{\la-\gamma}), \quad
\cL_{\la}^{-}(\zeta) = (1 \otimes \pi_{\zeta})(\cR_{\la-\gamma}^{(21)}).
$$
Using again the fact that $U_{\hbar}\N_{\pm}(\tau)^{\ge 2}$ is mapped to
zero by $\pi_{\zeta}$, we compute 
\begin{align} \label{chaat}
\nonumber L_{\la}^{+}(\zeta) & = 
\left( 1 +\theta(\hbar) f^{+}_{\la-\gamma
h+\gamma}(\zeta) \otimes E_{1,-1} \right) 
\left( k^{+}(\zeta+\gamma) 
\otimes E_{1,1}+ k^{+}(\zeta)^{-1} 
\otimes E_{-1,-1}
\right) 
\\ \nonumber
&  \left( 1 + \hbar e^{+}_{-\la}(z) \otimes E_{-1,1}\right) 
\\ & 
= 
\pmatrix 1 & \theta(\hbar) f^{+}_{\la - \gamma h +  \gamma}(\zeta)
\\ 0 & 1 \endpmatrix
\pmatrix k^{+}(z+\gamma) & 0 \\ 0 & k^{+}(z)^{-1} 
\endpmatrix
\pmatrix 1 & 0 \\ \hbar e^{+}_{-\la}(\zeta) & 1 \endpmatrix ,  
\end{align}
and
\begin{align*}
\cL^{-}_{\la}(\zeta) & \nonumber  = 
\left( 1 + \hbar e^{-}_{-\la}(\zeta) \otimes E_{-1,1}\right) 
q^{K\pa_{\zeta}}
\left( k^{-}(\zeta - \hbar) \otimes E_{1,1} + k^{-}(\zeta)^{-1} \otimes
E_{-1,-1} \right)
\\ &  
\left( 1 + \theta(\hbar) f^{-}_{\la - \gamma h +\gamma}(\zeta) \otimes
E_{1,-1}\right) 
\\ & \nonumber 
=
q^{K\pa_{\zeta}} L^{-}_{\la}(\zeta) 
, 
\end{align*}
where 
\begin{equation} \label{neila}
L^{-}_{\la}(\zeta) = \pmatrix 1 & 0 \\ \hbar e^{-}_{-\la}(\zeta-
K\hbar) & 1
\endpmatrix 
\pmatrix k^{-}(\zeta -\hbar) & 0 \\ 0 & k^{-}(\zeta)^{-1} \endpmatrix 
\pmatrix 1 & \theta(\hbar) f^{-}_{\la - \gamma h +\gamma}(\zeta ) 
\\ 0 & 1 \endpmatrix . 
\end{equation}

\begin{thm} \label{thm:RLL}
The matrices $L^{\pm}_{\la}(\zeta)$ defined by (\ref{chaat}) and
(\ref{neila}) satisfy the relations 
\begin{equation} \label{Lpm:Lpm}
R^{\pm}(\zeta-\zeta',\la) L^{\pm(1)}_{\la-\gamma h^{(2)}} (\zeta)
L^{\pm(2)}_{\la} (\zeta')
=
L^{\pm(2)}_{\la-\gamma h^{(1)}} (\zeta') L^{\pm(1)}_{\la} (\zeta)
R^{\pm}(\zeta-\zeta',\la-\gamma h) 
\end{equation}
\begin{align} \label{L+:L-}
L^{-(1)}_{\la}(\zeta) & R^{-}(\zeta-\zeta',\la-\gamma h)
L^{+(2)}_{\la}(\zeta')
\\ & \nonumber =
L^{+(2)}_{\la-\gamma h^{(1)}} (\zeta') R^{-}(\zeta-\zeta'+K\gamma,\la)
L^{-(1)}_{\la-\gamma h^{(2)}}(\zeta) {{A(\zeta,\zeta'-K\gamma)}
\over{A(\zeta,\zeta')}}. 
\end{align}
\end{thm}

{\em Proof. \/} It suffices to apply $id \otimes \pi_{\zeta} \otimes
\pi_{\zeta'}$, $\pi_{\zeta} \otimes \pi_{\zeta'} \otimes id$ and
$\pi_{\zeta} \otimes id \otimes \pi_{\zeta'}$ to (\ref{DYBE:ell}), after
the change of $\la$ into $\la-\gamma$, 
to simplify the coefficient $A(\zeta,\zeta')$ of Lemma \ref{adar} (which
is independent of $\la$), and to transfer the factors $q^{K\pa_{\zeta}}$
and $q^{K\pa_{\zeta'}}$ to the left. 
\hfill \qed \medskip 

\begin{remark} Connection of $A$ with $f_{K}$. The function
$A(\zeta,\zeta')$ of Lemma \ref{adar} satisfies the functional
equation 
$$
A(\zeta,\zeta')A(\zeta+\hbar,\zeta') =
{{\theta(\zeta-\zeta')}\over{\theta(\zeta-\zeta'+\hbar)}}; 
$$
after analytical prolongation, we see that $A$ only depends on
$\zeta-\zeta'$. The ratio $A(\zeta,\zeta' - K \gamma)
/ A(\zeta,\zeta')$ of relation (\ref{L+:L-}) is then simply connected
with the function $f_{K}$ expressing the commutator
$(k^{+}(z),k^{-}(w))$ (relation (\ref{k-k})) by 
$$
{{A(\zeta,\zeta'-K\gamma)}\over{A(\zeta,\zeta')}} =
f_{K}(\zeta'-\zeta-\hbar)^{-1}. 
$$
\end{remark}

\section{Elliptic quantum group $E_{\tau,\eta}(\frak{sl}_{2})$.}

\subsection{Definition}
 
Let us set $\eta = \gamma /2$ and define $E_{\tau,\eta}(\SL_{2})$
as the algebra generated by $h$ and the
$a_{i}(\la),b_{i}(\la),c_{i}(\la),d_{i}(\la),i\ge 0,\la\in\CC-L$, 
subject to the relations 
$$
[h,a_{i}(\la)] = [h,d_{i}(\la)] = 0, \quad
[h,b_{i}(\la)] = -2b_{i}(\la),  \quad [h,c_{i}(\la)] = 2c_{i}(\la),  
$$
and if we set 
$$
a(z,\la) = \sum_{i\ge 0} a_{i}(\la) e_{i,-\gamma h/2}(z), 
\quad b(z,\la) = \sum_{i\ge 0} b_{i}(\la) e_{i,\la-\gamma (h-2)/2}(z), 
$$
$$
c(z,\la) = \sum_{i\ge 0} c_{i}(\la) e_{i,-\la+\gamma (h+2)/2}(z),  
\quad d(z,\la) = \sum_{i\ge 0} d_{i}(\la)e_{i,\gamma h/2}(z),  
$$
and 
\begin{equation} \label{L-oper}
L(z,\la) = \pmatrix a(z,\la) & b(z,\la) \\ c(z,\la) & d(z,\la) 
\endpmatrix , 
\end{equation}
the relations
\begin{align} \label{RLL}  
 R^{+(12)}(z_{1}-z_{2}, \la - \gamma h) &  
 L^{(1)}(z_{1} , \la)  L^{(2)} (z_{2},\la - \gamma h^{(1)})  
\\ & \nonumber =  L^{(2)}(z_{2}, \la)  L^{(1)}(z_{1},\la-\gamma h^{(2)})  
 R^{+(12)} (\la, z_{1}-z_{2})   
\end{align} 
and 
\begin{equation} \label{Det=1}
\Det(z,\la) =
d(z+\gamma,\lambda)a(z,\lambda+\gamma)
 -
b(z+\gamma,\lambda)c(z,\lambda+\gamma)
{\theta(\lambda-\gamma h-\gamma)\over
\theta(\lambda-\gamma h)} = 1. 
\end{equation}

Here $R^{+}(z,\la)\in \End(\CC^{2} \otimes \CC^{2})$ is given by
(\ref{R+}); 
we define $h^{(1)}$ as $(E_{11} - E_{-1,-1}) 
\otimes 1$ and $h^{(2)}$ as $ 1 \otimes  
(E_{11} - E_{-1,-1})$. 
We also define as before, $f(\la-\gamma h)$ as $\sum_{\al\ge 0}
(\pa^{\al}_{\la}f)(\la)$ ${ {(-\gamma h)^{\al}} \over {\al!}}$,
 and 
$f(\la-\gamma h^{(i)})$ as $\sum_{\al\ge 0}
(\pa^{\al}_{\la}f)(\la)(-\gamma h^{(i)})^{\al}/\al!$. 

\begin{remark} The $L$-operator defined by (\ref{L-oper})
is $1$-periodic in the variables $z$ and $\la$, and satisfies 
\begin{equation} \label{per-L}
L(z+\tau,\la) = t_{\la-\gamma h} L(z,\tau) t_{\la}^{-1}, 
\end{equation}
where $t_{\la} = \pmatrix e^{-i\pi \la} & 0 \\ 0 & e^{i\pi \la}\endpmatrix$. 
On the other hand, the $R$-matrix (\ref{R-mat}) satisfies the conditions
\begin{equation} \label{per-R-2}
R^{+}(z+\tau,\la) = t^{(1)}_{\la-\gamma h^{(2)}} 
R^{+}(z,\la) t_{\la}^{-1(1)}. 
\end{equation}

The fact that the periodicity conditions (\ref{per-L}) and
(\ref{per-R-2})
seem compatible leads us to conjecture that the algebra $E_{\tau,\eta}
(\SL_{2})$ is a flat deformation of the function algebra of the group of 
holomorphic maps $L^{cl}$ from 
$(\CC - L)^{2}$ to $SL_{2}(\CC)$, such 
that  
$$
L^{cl}(z+\tau,\la) = \Ad(t_{\la})(L^{cl}(z,\la)). 
$$
Since the morphism $\Psi$ defined in Thm. \ref{comparison} is obviously 
surjective, this algebra has ``at least the size'' of 
$U_{\hbar}\G_{\cO}(\tau)$.  
\hfill\qed
\end{remark}

\subsection{Connection with the usual formulas}

The formulas defining the elliptic quantum groups in \cite{FV}
involve an $R$-matrix different from (\ref{R-mat}). Let us explain 
their connection with the above formalism. 

Consider the ring $\cF[[\gamma]]$ of formal series in $\gamma$, with  
coefficients meromorphic functions in $\la$. Let us adjoin to it a square  
root $\theta^{1/2}(\la)$ of $\theta(\la)$. The new ring $\cF_{1/2}[[\gamma]]$  
then contains a solution $\varphi$  
of the functional equation  
$$  
{{\varphi(\la + \gamma)}\over{\varphi(\la-\gamma)}} =  
{{\theta(\la)} \over{\theta(\la-\gamma)}};  
$$ 
we have  
$$ 
\varphi(\la) = \theta^{1/2}(\la)
\exp \left( - {1 \over {2\partial}}
\tanh (\gamma\partial/2) {{\theta'}\over{\theta}}
\right) (\la).  
$$ 

Define for $u,v \in \CC$, the quantities   
${\varphi(\la + \gamma(u h + v))}$ as  
$$ 
{\varphi(\la)}\exp\left( \sum_{k\ge 1} 
{{(\gamma(u h + v))^{k}}\over{k!}}  
\left({{\varphi'}\over{\varphi}}\right)^{(k-1)} 
\right) (\la). 
$$  
Note that any of the ratios ${{\varphi(\la + \gamma(u h + v))} 
\over {\varphi(\la + \gamma(u' h + v'))}}$, for $u,v,u',v' \in\CC$,  
belong to $\cF[h][[\gamma]]$.

\begin{lemma} \label{transfo} 
Let us set for $\la\in\CC - L$,  
\begin{equation} \label{a-abar-etc} 
\bar a(z,\la) =   {{\varphi(\la - \gamma h)}  
\over{\varphi(\la - \gamma )}}
a(z,\la),
\quad 
\bar b(z,\la) =   {{\varphi(\la - \gamma h)}  
\over{\varphi(\la + \gamma)}}b(z,\la),
\end{equation} 
\begin{equation} \label{c-cbar-etc} 
\bar c(z,\la) = {{\varphi(\la - \gamma h)}  
\over{\varphi(\la - \gamma)}} c(z,\la),
\quad 
\bar d(z,\la) =  {{\varphi(\la - \gamma h)}  
\over{\varphi(\la + \gamma)}} d(z,\la);
\end{equation} 
let us set $\bar L(z,\la) = \pmatrix \bar a(z,\la) & \bar b(z,\la) \\ 
\bar c(z,\la) & \bar d(z,\la)  
\endpmatrix$. Define 
\begin{align} \label{ell-R}
\bar R(z,\la) & =  E_{11} \otimes E_{11} + E_{-1,-1} \otimes E_{-1,-1}
+
{{\theta(\la+\gamma)\theta(z)} \over {\theta(\la)\theta(z-\gamma)}} 
E_{1,1} \otimes E_{-1,-1}
\\ & \nonumber
+ {{\theta(\la-\gamma)\theta(z)} \over {\theta(\la)\theta(z-\gamma)}} 
E_{-1,-1} \otimes E_{11}
-
{{\theta(\la+z)\theta(\gamma)} \over {\theta(\la)\theta(z-\gamma)}} 
E_{1,-1} \otimes E_{-1,1}
\\ & \nonumber 
- {{\theta(-\la+z)\theta(\gamma)} \over {\theta(-\la)\theta(z-\gamma)}} 
E_{-1,1} \otimes E_{1,-1} ; 
\end{align}
then we have the relations (see \cite{FV}) 
\begin{equation} \label{h-norm} 
h \bar a(z, \la) =  \bar a(z,\la) h, \quad  
h  \bar d(z,\la) =  \bar d(z,\la) h,  
\end{equation}  
\begin{equation} \label{h-norm2} 
h  \bar b(z, \la) = \bar b(z,\la) (h-2),  
\quad h \bar c(z,\la) =  \bar c(z,\la) (h+2),  
\end{equation}  
\begin{align} \label{bar-RLL}  
 \bar R^{(12)}(z_{1}-z_{2}, \la - \gamma h) &  
 \bar L^{(1)}(z_{1} , \la)  \bar L^{(2)} (z_{2},\la - \gamma h^{(1)})  
\\ & \nonumber =  \bar L^{(2)}(z_{2}, \la)  \bar 
L^{(1)}(z_{1},\la-\gamma h^{(2)})  
 \bar R^{(12)} (\la, z_{1}-z_{2}).   
\end{align} 
\end{lemma} 
 
{\em Proof.\/} 
We have  
$$ 
\bar R(z,\la)  = 
\varphi(\la - \gamma h^{(2)}) R^{+}(z,\la) \varphi(\la-\gamma h^{(1)})^{-1},  
$$ 
and 
$$ \bar L(z,\la)  = \varphi(\la - \gamma h)
L(z,\la) \varphi(\la-\gamma h^{(1)})^{-1} ;  
$$ 
Substitute these expressions in (\ref{RLL}); simplifications show 
that the $\bar R(z,\la)$ and $\bar L(z,\la)$ 
satisfy (\ref{bar-RLL}).  
\hfill \qed  \medskip 

\begin{remark}
The formulas of Lemma \ref{transfo} 
only use functions of $\cF[[\gamma]]$, although their proof 
uses the extension to $\cF_{1/2}[[\gamma]]$. 
\hfill \qed
\end{remark}
\medskip

\begin{remark} In \cite{FV}, the determinant is defined by
the formula
$$
\Det(z,\lambda)=\frac{\theta(\lambda)}{\theta(\lambda-\gamma h)}
(\bar d(z+\gamma,\lambda)\bar a(z,\lambda+\gamma)
-
\bar b(z+\gamma,\lambda)\bar c(z,\lambda+\gamma)).
$$
This formula is equivalent to the first equation of (\ref{Det=1}), 
as one sees by inserting the expressions
for $\bar a,\dots,\bar d$ in terms of $a,\dots,d$ and using the
identity
$$
\frac{\varphi(\lambda-\gamma h)\varphi(\lambda-\gamma h+\gamma)}
{\varphi(\lambda)\varphi(\lambda+\gamma)}
=\frac{\theta(\lambda-\gamma h)}{\theta(\lambda)}\, .
$$
\end{remark}

\begin{remark}
By tensoring them with $1$-dimensional representations, we can view the 
evaluation representations studied in \cite{FV} as representations of
the factor algebra introduced in this paper by the relation
$\Det(z,\la)=1$. After expansion in series in $\eta = \gamma /2$,
the formulas defining the evaluation representations of \cite{FV} only
have singularities for $\la\in L$ or $z-w\in L$. The effect of the
tensoring with $1$-dimensional representations is to multiply the matrix
$L(w,\la)$ by a function $g_{z}(w)$ satisfying
$$
g_{z}(w+\gamma)g_{z}(w)=
{\theta(z-w+(\Lambda-1)\eta)\over
\theta(z-w-(\Lambda+1)\eta)}, \quad \eta=\gamma/2; 
$$
this equation can be solved in a similar way to that for $\varphi$, and
we will find for $g_{z}(w)$ a formal series in $\gamma$ with
coefficients functions of $z-w$ with only singularities for
$z-w\in L$. 

Therefore the final representations can be viewed as
representations of the algebras $E_{\tau,\eta}(\SL_{2})$, provided $w$
is considered as a formal variable at the origin 
(as it is the case in \cite{HGQG}). 
\end{remark}

\section{Quantum currents for $E_{\tau,\eta}(\SL_{2})$}

\begin{thm} \label{comparison}
There is a morphism $\Psi$ from $E_{\tau,\eta}(\SL_{2})$ to 
$U_{\hbar}\G_{\cO}(\tau)$, defined by the formulas
$$
\Psi(h) = h , 
$$
\begin{equation} \label{bar-a}
\Psi (  a(z, \la) ) = \hbar\theta(\hbar)
f^{+}_{\la - \hbar - \gamma h}(z) k^{+}(z)^{-1}
e^{+}_{-\la}(z)
+ k^{+}(z - \hbar )
\end{equation}
\begin{equation}  \label{bar-b}
\Psi ( b(z,\la) )
= \theta(\hbar) f^{+}_{\la - \hbar -\gamma h}(z) k^{+}(z)^{-1} 
\end{equation}
\begin{equation} \label{bar-d-c}
\Psi (  d(z, \la) )  = k^{+}(z)^{-1}, \quad  \Psi( c(z,\la) ) = \hbar
k^{+}(z)^{-1} e^{+}_{-\la}(z),
\end{equation}
\end{thm}

{\em Proof. \/} Let us first show that these formulas are generating series
for images of the $a_{i},b_{i},c_{i},d_{i},i\ge 0$. For this, we note that 
their right-hand sides are holomorphic functions on $(\CC - L)^{2}$ with 
values in $U_{\hbar}\G_{\cO}(\tau)$, $1$-periodic in $z$ and $\la$ and with 
the quasi-periodicity
properties in $z$ and $\la$ described by (\ref{per-L}). For the periodicity 
properties in $z$, this follows directly from (\ref{x-per}), the equations 
$$
k^{+}(z+1) = k^{+}(z), \quad k^{+}(z+\tau) = e^{i\pi\gamma h} k^{+}(z),  
$$
which are proved in the same way as (\ref{K-per}), and the commutation 
relations between $h$ and the $e^{+}_{\la}(z), f^{+}_{\la}(z)$. 

By Thm. \ref{thm:RLL}, $\Psi(a(z,\la)), 
\Psi(b(z,\la)), \Psi(c(z,\la))$ and $\Psi(d(z,\la))$ satisfy the 
relations (\ref{RLL}). Finally, one computes that the image by $\Psi$ of the 
middle term of equation (\ref{Det=1}) is equal to $1$. 
This ends the proof of the theorem. 
\hfill \qed\medskip

\begin{remark} Since $\Psi(d(z,\la))$ is independent of 
$\la$, it should be clear that $\Psi$ is not
injective. 
\hfill \qed \medskip 
\end{remark}

\begin{remark}
There is an algebra morphism from the tensor product
$E_{\tau,\eta}(\SL_{2})\otimes \Diff$ $(\CC - L)$ to 
$E_{\tau,\eta}(\SL_{2})^{\otimes 2}\otimes \Diff(\CC - L)$; 
the formulas for it are $\Delta(L(z,\la)) = L^{(13)}(z-\gamma h^{(2)}, \la)
L^{(23)}(z, \la)$. It would be interesting to understand better the
relation of this formula with (\ref{Del12R}).  
\hfill \qed \medskip 
\end{remark}

\begin{remark}
Relations (\ref{Lpm:Lpm}) and  (\ref{L+:L-})
suggest to define double elliptic quantum groups
generated by the matrices $L^{\pm}(z,\la)$, the derivation $D$ and the
central element $K$, with the following functional properties: 
the $L^{\pm}(z,\la)$ are holomorphic functions in the variable 
$\la\in\CC - L$; $L^{+}(z,\la)$ also depends holomorphically on $z\in
\CC - L$, and $L^{-}(z,\la)$ is a a regular
power series in $z$; the periodicity conditions for
$L^{\pm}(z,\la)$ in $\la$ are the same as those for
$L(z,\la)$, and the periodicity conditions for $L^{+}(z,\la)$ in $z$ are
the same as those for $L(z)$; 
and satisfying relations (\ref{Lpm:Lpm}), (\ref{L+:L-})
and $[D,L^{\pm}(z,\la)] = \pa L^{\pm}(z,\la)/\pa z$. This
algebra should be, as it is the case for
$E_{\tau,\eta}(\SL_{2})$ with respect to $U_{\hbar}\G_{\cO}(\tau)$,
somewhat larger than $U_{\hbar}\G(\tau)$.
\end{remark}

\end{document}